%% file: integral.tex
\documentclass{llncs}

\usepackage[T1]{fontenc}

\usepackage{graphicx}

\usepackage{amsfonts,amssymb,amsmath}

\usepackage[linesnumbered,vlined,boxed,fillcomment]{algorithm2e}
\SetKw{In}{ in }
\SetKw{And}{ and }
\SetKw{All}{ all }
\SetKw{Or}{ or }
\SetKw{Break}{break}
\SetKw{Continue}{continue}
\SetKw{KwFrom}{from }
\makeatletter
\let\oldnl\nl% Store \nl in \old >nl
\newcommand{\nonl}{\renewcommand{\nl}{\let\nl\oldnl}}% Remove line number for one line
\makeatother

\usepackage{enumerate}%
\usepackage{tikz}%
\usetikzlibrary{calc,positioning}
\usetikzlibrary{decorations.markings}%

\usepackage{enumitem}

\usepackage{color}
\usepackage{etoolbox}
\makeatletter
\let\llncs@addcontentsline\addcontentsline
\patchcmd{\maketitle}{\addcontentsline}{\llncs@addcontentsline}{}{}
\patchcmd{\maketitle}{\addcontentsline}{\llncs@addcontentsline}{}{}
\patchcmd{\maketitle}{\addcontentsline}{\llncs@addcontentsline}{}{}
\makeatother
\usepackage{hyperref}%

\hypersetup{% gestion des liens
  colorlinks=true,% colorise les liens
  breaklinks=true,% permet le retour a la ligne dans les liens trop longs
  urlcolor= blue,% couleur des hyperliens
  linkcolor= blue,% couleur des liens internes
  bookmarksopen=true,% si les signets Acrobat sont crees, les afficher completement.
  pdfpagemode=UseOutlines,% afficher la table des matieres dans la navigation du pdf ?
  bookmarksopenlevel=1,%
  pdftitle={Integral Bases},% informations apparaissant dans
  pdfauthor={Adrien Poteaux \& Martin Weimann},% dans les informations du document
}

\input{commands.tex}

\begin{document}

\title{Fast Integral Bases Computation}

% les orcid apparaissent en clair dans le papier, c'est moche, je commente en attendant
\author{Adrien Poteaux\inst{1}%\orcidID{0000-0002-7493-3001}
  \and Martin Weimann\inst{2}}%\orcidID{0000-0002-8109-5659}
\institute{Univ. Lille, CNRS, Centrale Lille, UMR 9189 CRIStAL, F-59000 Lille, France\\%
  \email{adrien.poteaux@univ-lille.fr}\\%
  \url{https://www.fil.univ-lille.fr/~poteaux/}
  \and   LMNO, Universit\'e de Caen-Normandie\\%
  \email{martin.weimann@unicaen.fr}\\%
  \url{https://weimann.users.lmno.cnrs.fr/}
}

\maketitle

\begin{abstract}
  We obtain new complexity bounds for computing a triangular integral
  basis of a number field or a function field. We reach for function
  fields a softly linear cost with respect to the size of the output
  when the residual characteristic is zero or big enough. Analogous
  results are obtained for integral basis of fractional ideals, key
  ingredients towards fast computation of Riemann-Roch spaces. The
  proof is based on the recent fast OM algorithm of the authors and on
  the \maxmin algorithm of Stainsby, together with optimal
  truncation bounds and a precise complexity analysis.
  \keywords{Fields, Fractional ideals, Integral bases, OM algorithm.}
\end{abstract}

\section{Introduction}
\label{sec:intro}
\input{intro}

\section{OM algorithm and local integral basis}
\label{sec:locbasis}
\input{localbasis}

\subsection{Precision of the OM algorithm}
\label{sec:precision}
\input{precision}

\section{Triangular $p$-integral bases}
\label{sec:p_int_basis}
\input{triangular}

\subsection{Global integral bases}
\label{sec:global_basis}
\input{global}

\section{Bases of fractional ideals}
\label{sec:fractional}
\input{fractional}

\section{Complexity and proofs of the main results}
\label{sec:complexity}
\input{complexity}
\section{An illustrative example}
\label{sec:example}
\input{stainsby-3.3}

\bibliographystyle{splncs04} {\bibliography{tout}}
%\addcontentsline{toc}{section}{References.}

\end{document}

%% file: commands.tex
\newcommand{\hbb}[1]{\ensuremath{\mathbb{#1}}}% anneaux / corps

\newcommand{\Fi}{\hbb{F}}

\newcommand{\Ni}{\hbb{N}}
\newcommand{\Qi}{\hbb{Q}}

\newcommand{\Zi}{\hbb{Z}}

\newcommand{\hcal}[1]{\ensuremath{\mathcal{#1}}}
\newcommand{\Ac}{\hcal{A}}
\newcommand{\Bc}{\hcal{B}}

\newcommand{\Fc}{\hcal{F}}

\newcommand{\Nc}{\hcal{N}}
\newcommand{\Oc}{\hcal{O}}

\newcommand{\Pc}{\hcal{P}}

\newcommand{\hfrak}[1]{\ensuremath{\mathfrak{#1}}}
\newcommand{\q}{\hfrak{q}}% Number of elements in the ground field
\newcommand{\p}{\hfrak{p}}

\newcommand{\hbf}[1]{\ensuremath{\mathbf{#1}}}

\newcommand{\ib}{\hbf{i}}
\newcommand{\jb}{\hbf{j}}

\newcommand{\eb}{\hbf{e}}

\newcommand{\maxmin}{\textnormal{\texttt{MaxMin}}}

\renewcommand{\O}{\Oc}% O
\newcommand{\Ot}{\O\widetilde{\ }}% soft-O
\newcommand{\Oe}{\O_\eps}% soft-O
\newcommand{\dx}{{d}}% degre en Y de F
\renewcommand{\char}{\textrm{char}}% caracteristique du corps

\newcommand\disc[1][x]{{\textrm{Disc}_{#1}}} % je pense qu'il ne faut pas utiliser \Delta pour ne pas s'embrouiller avec les aretes du polygone de Newton
\newcommand{\Card}{\textrm{Card}}
\newcommand{\Res}{\textrm{Res}}

\newcommand{\ind}{\textrm{ind}}

\newcommand{\war}[1][w]{\stackrel{\rightarrow}#1}%w arrow

\newcommand{\eps}{\epsilon}% pour changer facilement \epsilon / \varepsilon

%% file: intro.tex
Let $A$ be a principal ideal domain with field of fractions $K$ and
let $L/K$ be a separable field extension of degree $d$. Denote $B$ the
integral closure of $A$ in $L$. It is a free $A$-module of rank $d$
and an integral basis of $L/K$ is a collection $b_1,\ldots,b_d\in B$
that form a $A$-basis of $B$.

This paper intends to give new complexity bounds for computing an
integral basis and, more generally, for computing an $A$-basis of an
arbitrary fractional ideal of $L$. These fundamental problems of
computer algebra are cornerstones towards more advanced computational
issues both in algebraic number theory and in algebraic geometry, such
as developing a fast arithmetic of ideals in number fields \cite{Co12}
or computing Riemann-Roch spaces in function fields \cite{Ba14,He02}.

In all of the sequel, we will assume that $L=K(\theta)$ is
  generated by the root $\theta$ of a degree $\dx$ monic separable
  irreducible polynomial $f\in A[x]$. Under this assumption, there
exists an integral basis of $L/K$ of shape
\begin{equation}\label{eq:baseB}
  \Bc=\left(1,
    \frac{g_1(\theta)}{a_1},\ldots,\frac{g_{d-1}(\theta)}{a_{d-1}}\right)
\end{equation}
where $g_i\in A[x]$ is monic of degree $i$, $a_i\in A\setminus\{0\}$
and $a_1\,|\, a_2 \,|\cdots |\, a_{d-1}$. We call such a basis a
\emph{triangular integral basis}. This specific shape has various
algorithmic advantages, in particular for computing Hermite normal
forms or Popov forms.

% \begin{enumerate}
% \item Compute an integral basis of $L$, that is a collection
%   $b_1,\ldots,b_d\in B$ that form a basis of $B$ as an $A$-module.
% \item Compute an $A$-basis of any given fractionnary ideal $I$ of
%   $L$ (in particular, the conductor ideal of $B$ in $A$).
% \end{enumerate}
%
% Notice that task one is a key point in the computation of
% Riemann-Roch spaces of a plane curve.

% Let $\theta\in B$ stands for a root of $F$. We have an inclusion
% $A[\theta]\subset B$ of free $A$-modules of rank $d$. Hence, the
% $A$-module $B/A[\theta]$ is torsion. Assuming that $A$ is PID, we
% may thus write
%\[
%B/A[\theta]\simeq \sum_{i=1}^d A/a_i A
%\]
%for some $a_i\subset A\setminus\{0\}$. The ideal $(a_1\cdots a_d)$
% only depends on $F$.

\paragraph{Model of computation.} For the sake of simplicity, we will
express our complexity results in the case $A=k[t]$
(resp. $A=\Zi$). We work with computation trees \cite[Section
4.4]{BuClSh97}. When considering $A=k[t]$, we use an algebraic RAM
model, counting only the number of arithmetic operations in $k$. For
$A=\Zi$, we consider the Boolean RAM model.

We classically denote $\O()$ and $\Ot()$ to respectively hide constant
and logarithmic factors in our complexity results ; see
e.g. \cite[Chapter 25, Section 7]{GaGe13}. We additionally let
$\Oe(\dx)=\O(d^{1+\eps(\dx)})$ with $\eps(\dx)\to 0$. We have
$\Ot(\dx)\subset \Oe(\dx)$, and freely speak of \emph{almost linear}
in $\dx$ for both notations. As in \cite{PoWe22b}, we express our
complexity with $\Oe()$ because we deal with dynamic evaluation via
the deterministic algorithm of \cite{HoLe19}. Our results could be
express with the $\Ot()$ notation using some Las Vegas sub-algorithm
instead.

\paragraph{Discriminant and index.}
Our results will be expressed using the size functions $h(a)=\deg(a)$
when $A=k[t]$ (resp. $h(a)=\log(|a|)$ when $A=\Zi$). Note that
$h(a)=\sum_{p} v_p(a) h(p)$ where the sum runs over the primes $p$ of
$A$ and $v_p$ stands for the $p$-adic valuation. Let us denote
$\Delta=\Delta_f$ the discriminant of $f$ and $\Delta_{L/K}$ the
discriminant of $L/K$. Both quantities are related by the formula
\begin{equation}\label{eq:disc}
  \Delta =D^2 \,\Delta_{L/K}
\end{equation}
where $D=D_f\in A$ is the index of $f$, generator of the index ideal
$[B:A[\theta]]$. We denote for short $h_\Delta:=h(\Delta)$,
$h_D:=h(D)$ and $h_{red}:=\sum_{p^2|\Delta} h(p)$.  Looking at the
triangular basis \eqref{eq:baseB}, we observe that
$D= u a_1\cdots a_{d-1}$ for some unit $u\in A^\times$ and, up to
reduce $g_i$ modulo $a_i$, the basis $\Bc$ has global size
$\Oc(d h_D)$, this bound being sharp.
% (although it's not clear if this bound remains sharp after Hermite
% reduction \adrien{J'ai un peu de mal avec cette phrase, peut \^etre
% insister que c'est sharp dans le cas triangulaire ?  Bref, \`a
% discuter.}).
This observation suggests to emphasise the dependency of the
complexity on the index $D$ instead of the discriminant $\Delta$ which
is classically considered in the literature (see
e.g. \cite{Abe20,Ba16,St18}).

\paragraph{Global integral basis computation.} Our main result is:

\begin{theorem}\label{thm:main_hD}
  Suppose given a squarefree factorisation of the index $D$. There
  exists a deterministic algorithm which computes a triangular
  integral basis of $L/K$ in less than
  \begin{enumerate}
  \item\label{it:main1} $\Oc_\eps(\dx h_D)$ operations in $k$ if
    $A=k[t]$ with $\char(k)=0$ or $\char(k)>\dx$,
  \item $\Oc_\eps(\dx\,h_D+h_D^2)$ operations in $k$ if $A=k[t]$ with
    $\char(k)\leq\dx$,
  \item $\Oc_\eps(\dx\,h_D+h_D^2)$ word operations if $A=\Zi$.
  \end{enumerate}
  If we are only given the squarefree factorisation of the
  discriminant $\Delta$, then a similar result holds, adding a cost
  $\Oc(\dx \,h_{red})$ to the complexity estimates.
\end{theorem}
  
We thus get a softly linear cost in case \ref{it:main1}. In the other
cases, the over-cost in $h_D^2$ is only due to the complexity of the OM
algorithm \cite{PoWe22b} over small residual characteristic, so that
\emph{the result would become softly linear for all cases if we are
  able to improve the complexity of the OM algorithm}. In the last
statement, the extra cost $\Oc(\dx \,h_{red})$ is needed to detect
which prime divisors of $\Delta$ are prime divisors of $D$.

In the case $A=k[t]$, we have access to polynomial time squarefree
factorisation, leading to the following total cost estimate:

\begin{corollary}
  \label{thm:main_hDelta} Suppose $A=k[t]$. There exists a
  deterministic algorithm which computes a triangular integral basis
  of $L/K$ in less than
  \begin{enumerate}
  \item $\Oc_\eps(\dx \,h_\Delta)$ operations in $k$ if $\char(k)=0$
    or $\char(k)>\dx$,
  \item $\Oc_\eps(\dx \,h_\Delta +h_D^2)$ operations in $k$ otherwise.
  \end{enumerate}
  The algorithm returns as a byproduct the index $D$ and the
  discriminant $\Delta_{L/K}$.
\end{corollary}

These results improve significantly \cite[Theorem 3.5]{St18} which
after Chinese remainder gluing leads to a global cost
$\Oc_\eps(\dx^2\,h_\Delta+\dx\,h_\Delta^2)$.  Up to our knowledge,
this is the best complexity estimate in the literature. Note also
\cite{Abe20} which leads to $\Oc_\eps(\dx^2\,h_\Delta)$ in the
particular case $A=k[t]$ and $\char(k)=0$, based on fast computation
of Puiseux series (although the resulting integral basis is not
triangular).

The complexity indicators $h_D$ and
$h_{red}$ in Theorem \ref{thm:main_hD} satisfy
$h_D + h_{red} \le h_\Delta$, and the difference may be
significant, especially when wild ramification occurs.

\begin{example}\label{xmp:tame}
  Consider $f=x^q+t^{n} x +t\in \Fi_q[t][x]$ with $q$ a
  prime. We have $\Delta=t^{n q}$ while $D=1$. We get $h_\Delta=q n$,
  $h_D=0$, and $h_{red}=1$.  An integral basis is trivially
  $(1,\theta,\ldots,\theta^{d-1})$. Accordingly, Theorem
  \ref{thm:main_hD} estimates the true cost $\Oc(d)$ to compute such a basis, while using
  $h_\Delta$ as a complexity indicator would lead to the bad estimate
  $\Oc_\eps(d^2 n)$, the integer $n$ being arbitrarily large. See \cite[Section 6]{GuMoNa13} for further examples
  comparing $h_D + h_{red}$ and $h_\Delta$ in the case $A=\Zi$.
\end{example}

Let us emphasise that Corollary \ref{thm:main_hDelta} simply adds to
Theorem \ref{thm:main_hD} an extra cost $\O(d h_\Delta)$ due to the
computation of $\Delta$, plus $\Ot(h_\Delta)$ for its squarefree
factorisation. Hence any progress for computing the discriminant would
improve the complexity estimate of Corollary
\ref{thm:main_hDelta}. There have been recent results in that
direction in the case $k=\Fi_q$ : it follows from \cite{Vil23} that
there is a randomised algorithm of Monte Carlo type which computes the
radical of $\Delta$ (and this is enough for our purpose) in softly
optimal time in $h_\Delta$.

\paragraph{Integral basis over a prime $p$.} If we only want a
$p$-integral basis at a given prime $p\in A$ (that is an $A_p$-basis
of $B\otimes A_p$ over the localisation of $A$ at $p$), we get similar
complexity estimates than in Theorem \ref{thm:main_hD}, replacing
$h_D$ and $h_{red}$ by their respective local contributions
$v_p(D)h(p)$ and $h(p)$. There's no need to compute and factorise the
discriminant in such a case. This result has to be compared to
\cite[Theorem 3.5]{St18} (triangular basis) or \cite[Lemma 3.10]{Ba16}
(non triangular basis).

\paragraph{The case of fractional ideals.} In Section
\ref{sec:fractional}, we provide similar complexity results  to compute integral basis of an arbitrary
fractional ideal $I\subset L$, expressed now in terms of the size of
the index $[I^*:A[\theta]]$ of the smallest multiple $I^*=\alpha I$, $\alpha\in K$ such that $B\subset I^*$, see Theorem  \ref{thm:main_fractional}. This is an important issue towards the computation of Riemann-Roch spaces in function fields.

\paragraph{Previous results.} Classical methods for computing integral bases are variants of the
Round-2 and Round-4 routines by Zassenhaus and Ford
\cite{FoPaRo10,Hal01,Ho94,Poh93}. The central idea is to start from a
known order $A[\alpha]\subset B$, and to enlarge it for each prime
$p\in A$ dividing the discriminant of $A[\alpha]$ until we reach the
maximal order $B$. Another strategy of local-to-global type was
developed by Okutsu \cite{Ok82}: assuming given the local
$\p$-integral basis $\Bc_\p$ for each prime ideals $\p\subset B$
dividing a prime $p\in A$, one can compute multipliers $z_\p\in L$
such that $\Bc_p=\cup z_\p \Bc_\p$ is a $p$-integral basis. Then,
after reducing $\Bc_p$ to a triangular form, one can glue the various
bases $\Bc_p$ into a global integral basis $\Bc$ by means of Chinese
remainders.  Later on, Montes \cite{Mo99} extended the ideas of Ore
and MacLane \cite{Ma36b,Ma36a} and developed the OM algorithm that
computes a representation of the prime ideals $\p$ dividing a prime
$p\in A$ by way of factoring $f$ over the $p$-adic completion of
$K$. This led to an efficient computation of the local
basis $\Bc_\p$, last missing ingredient in Okutsu's strategy. This
opened the door to various OM-based routines: methods of multipliers
\cite{Abe20,Ba16,GuMoNa13}, method of the quotients \cite{GuMoNa15}
and finally the remarkably simple \maxmin{} algorithm of Stainsby
\cite{St18} that we follow here.
% which allows to glue $\p$-basis in a \emph{triangular} $p$-integral
% basis avoiding any costly Hermite type reduction.  Coupled with the
% recent fast divide and conquer OM-algorithm \cite{PoWe22b}, this
% leads to the new complexity estimates presented in this paper.

\paragraph{Summary of our strategy.} We do not claim any originality
in our approach, following the classical local to global strategy
outlined in \cite{Ba14,St18}:
\begin{enumerate}
\item Run the fast OM algorithm \cite{PoWe22b} above each prime
  $p\in A$ dividing $\Delta$ (if we start with a squarefree factor,
  we rely on dynamic evaluation) and deduce for each prime ideal
  $\p\subset B$ dividing $p$ a local $\p$-integral basis together with
  a suitable approximant of the associated local factor $F_\p$ of $f$.
\item \label{it:strat2} Apply the \maxmin{} algorithm \cite{St18} to
  deduce a triangular $p$-integral basis (each numerator
  being a suitable multiplicative combination of the numerators of the
  various $\p$-bases).
\item Use the Chinese Reminder Theorem to glue these $p$-integral
  basis as a global triangular integral basis.
\end{enumerate}
Concerning Point \ref{it:strat2}, it is remarkable that \cite{St18}
allows to compute a \emph{triangular} $p$-integral basis avoiding the
usual costly Hermite type reduction step. The OM algorithm was the
main bottleneck of the associated complexity analysis given in
\cite{St14} and our improvements mainly follow from the recent faster
OM algorithm \cite{PoWe22b}, together with a careful study of the
various $p$-adic precisions needed to conduct the computations.

\paragraph{Organisation of the paper.}
In Section \ref{sec:locbasis}, we first remind how local integral
bases (one for each factor of the local factorisation of $f$) are
deduced from the Okutsu frames computed during the OM algorithm. Then,
we prove that the usual precision $v_p(\Delta)$ used in the OM
algorithm can be improved by $2\, v_p(D)+1$, key result towards the
proof of Theorem \ref{thm:main_hD}. Section \ref{sec:p_int_basis} is
dedicated to glue local bases into a reduced triangular $p$-integral
basis thanks to the \maxmin{} algorithm of Stainsby, these $p$-integral bases being then glued into
an integral basis of $L/K$ using Chinese remaindering (Subsection \ref{sec:global_basis}). Section
\ref{sec:fractional} is dedicated to integral bases of arbitrary
fractional ideals. We pay attention to complexity issues and prove our
main Theorem \ref{thm:main_hD} in Section \ref{sec:complexity}
%, providing also some computational examples
. Finally, we illustrate the
overall strategy with an example of \cite{St18} in Section
\ref{sec:example}.

%%% Local Variables:
%%% mode: latex
%%% TeX-master: "integral"
%%% End:

%% file: localbasis.tex
Let $A$ be a principal ideal domain with field of fraction $K$ and let
$L=K(\theta)$ be the field extension determined by a root
$\theta\in \overline{K}$ of a monic, irreducible and separable
polynomial $f\in A[x]$ of degree $d$.

The OM algorithm \cite{Mo99}, from the initial of its main artisans
Ore, MacLane, Okutsu and Montes, computes an OM representation of each
prime ideal of $L$ lying over a given prime $p\in A$. This
computational object supports several arithmetic data attached to an
irreducible factor (say $F$) of $f$ over the $p$-adic completion $K_p$
of $K$. In particular, it allows to compute a local integral basis of
the finite extension of $K_p$ determined by $F$, first step towards
the local-to-global computation of an integral basis of $L/K$.  After
recalling this construction, we give new tight bounds for the
precision required by the OM
algorithm. %, key ingredient towards the proof of Theorem \ref{thm:main_hD}.

\subsection{Local integral basis and OM-factorisation.}
We fix a prime $p\in A$ and consider $v_p:K^\times \to \Zi $ the
$p$-adic valuation, $K_p$ the completion of $K$ with respect to $v_p$
and $\Oc_p$ the valuation ring of $K_p$. We still denote $v_p$ the
canonical extension of $v_p$ to a fixed algebraic closure
$\overline{K}_p$ of $K_p$.

\paragraph{Okutsu frame and local integral basis.} Let $F\in \Oc_p[x]$
be an irreducible monic polynomial of degree $n$ and let
$\alpha\in \overline{K}_p$ be a root of $F$.

\begin{definition}\label{def:okutsu}
  An \emph{Okutsu frame} of $F$ is a sequence
  $[\phi_1,\ldots,\phi_{r+1}]$, with $\phi_i\in\Oc_p[x]$ monic of
  degree $m_i$ such that (denoting $m_0=1$, $\phi_0=1$):
  \begin{itemize}[label=$\bullet$]
  \item $m_1 \,|\, m_2 \,|\,\cdots \,|\, m_{r+1}$ and
    $1\le m_1<m_2 < \cdots <m_{r+1}=n$.
  \item For all $g\in \Oc_p[x]$ monic,
    $\deg g < m_{i+1} \Longrightarrow \frac{v_p(g(\alpha))}{\deg(g)}
    \le \frac{v_p(\phi_i(\alpha))}{m_i}<
    \frac{v_p(\phi_{i+1}(\alpha))}{m_{i+1}}$
  \end{itemize}
  The polynomial $\phi_{r+1}$ is called an \emph{Okutsu approximant}
  of $F$.
\end{definition}

% Okutsu frames are crucial objects as they capture the main
% invariants of an irreducible polynomial. In particular,

The degrees $m_i=\deg(\phi_i)$ and the length $r$ do not depend on the
choice of the frame. For any $0\le m <n$, we can write in a unique way
$m=j_0 \,m_0+\cdots + j_r \,m_r$ with $0\le j_i< m_{i+1}/m_{i}$. We
accordingly let
% \begin{equation}\label{eq:divisor_polynomials}
\[
  g_m:=x^{j_0}\phi_1^{j_1}\cdots \phi_r^{j_r}\in \Oc_p[x].
\]
%\end{equation}

The integral closure $\overline{\Oc}_p$ of $\Oc_p$ in $K_p(\alpha)$ is
a free $\Oc_p$-module of rank $n$. Okutsu proved \cite[Theorem
1]{Ok82}:
\begin{proposition}\label{prop:local_basis}
  Let $\eta_m=\lfloor v_p(g_m(\alpha)) \rfloor$. The family
  $ 1\, ,\, \frac{g_1(\alpha)}{p^{\eta_1}}\, ,\, \ldots\,
  ,\,\frac{g_{d-1}(\alpha)}{p^{\eta_{d-1}}} $ is an $\Oc_p$-basis of
  $\overline{\Oc}_p$. We call it an \emph{Okutsu basis} of $F$.
\end{proposition}
We have $\eta_1 \le \cdots \le \eta_{d-1}=\exp(F)$ where $\exp(F)$ is
the \emph{integrality exponent} of $F$, least integer such that
$p^{\exp(F)}\overline{\Oc}_p \subset \Oc_p[\alpha]$.

\begin{definition}\label{def:okutsu_num}
  Let $g_0=1$ and $g_n=\phi_{r+1}$. The set
  $\Nc(F):=\{g_0,g_1,\ldots,g_{n-1},g_n\}$ is called an extended set
  of \emph{Okutsu numerators} of $F$.
\end{definition}

\paragraph{OM-factorisation.}\label{ssec:OM_facto}

The prime ideals $\p$ of $L$ dividing $p$ are one-to-one with the
irreducible monic factors of $f$ in $\Oc_p[x]$. We denote
$f=\prod_{\p|p} F_\p$ and we let $\theta_\p\in \overline{K}_p$ be an
arbitrary root of $F_\p$.

\begin{definition} 
  An Okutsu factorisation of $f$ above $p$ is a set $(\Fc_\p)_{\p|p}$
  where $\Fc_\p:=\left[\phi_{\p,0},\ldots,\phi_{\p,r_{\p}+1} \right]$
  is an Okutsu frame of $F_\p$, with $\phi_{\p,i}\in
  A[x]$. An OM-factorisation of $f$ above $p$ is an Okutsu
  factorisation s.t. the approximants
  $\phi_\p:=\phi_{\p,r_{\p+1}}$ satisfy
  $v_p(\phi_\p(\theta_\p))>v_p(\phi_\p(\theta_\q))$ for all $\p\ne \q$
  (see \cite[Definition 3.2]{BaNaSt13}).
\end{definition}

The stronger condition for being an OM-factorisation ensures that the
approximant $\phi_\p$ uniquely determines the corresponding factor
$F_\p$ of $f$ (we might have $\phi_\p=\phi_\q$ for some $\q\ne \p$ in
an Okutsu factorisation).

The fast OM algorithm of \cite{PoWe22b} will allow us to compute an
OM-factorisation of $f$ in the aimed complexity bound, together with
all the $\phi_{\p,i}$ defined above, thus an Okutsu basis of each
$F_\p$ thanks to Proposition \ref{prop:local_basis}. 
%Assuming that the
%precision of the approximants is high enough, this will enable the
%computation of a triangular $p$-basis in Section
%\ref{sec:p_int_basis}.
% We now provide a tight bound for this required precision, namely
% $2\,v_p(D)+1$, improving the usual $v_p(\Delta)$.

%%% Local Variables:
%%% mode: latex
%%% TeX-master: "integral"
%%% End:

%% file: precision.tex
We improve here the results of \cite{BaNaSt13} about the precision
required for computing an OM-factorisation of $f$ (and a $p$-integral
basis). This section is quite technical, depending a lot of the
references \cite{BaNaSt13,Na14,PoWe22b}. It is independent of the
strategy for the computation of a triangular integral basis and the
reader mainly interested on that part can skip it.

\paragraph{The Okutsu bound.}  For $F\in \Oc_p[x]$ monic, separable
and irreducible, and $\alpha$ a root of $F$, we define the
\textit{Okutsu bound} of $F$ as
%\begin{equation}\label{eq:okutsu_bound}
\[
  \delta_0(F): = \deg(F) \;
  \max \left\{\frac{v_p(g(\alpha))}{\deg(g)},\; g\in \Oc_p[x]\,\,\mathrm{monic},\, \deg(g)<d\right\}.
\]
%\end{equation}
Given $f\in A[x]$ with irreducible factorisation
$f=\prod_{\p|p} F_\p\in \Oc_p[x]$ as above, and denoting
$d_\p=\deg(F_\p)$, we let
\[
  \delta^*(f):=\frac{1}{2}\sum_{\p|p} d_\p \delta_0(F_\p) +
  \sum_{\p\ne \q} v_p(Res(F_\p,F_\q)).
\]
In contrast to the discriminant valuation $\delta(f):=v_p(\Delta(f))$,
the rational number $\delta^*(f)$ is an Okutsu invariant of $f$, that
is it only depends on combinatorial data attached to an
OM-factorisation $f$ (see \cite{Na14} for details). Both quantities
are related by Proposition \ref{prop:deltastar} below.

\paragraph{Bounds for $\delta^*(f)$.} For each $\p | p$, let $e_\p$
and $f_\p$ be respectively the ramification index and the residual
degree of $\p$ over $p$ and let $L_\p=K_p(\theta_\p)$. It is well known that $\Delta_{L_\p/K_\p})\ge f_\p(e_\p-1)$, with equality if and only if $L_\p/K_p$ is
tame, that is if $p\nmid e_\p$ and the 
residue field extension is separable. We define
\[
  \rho(f):=\sum_{\p|p}
  \left(v_p(\Delta_{L_\p/K_\p})-f_\p(e_\p-1)\right)\ge 0
\]
which thus measures the non tameness of $L/K$.
%We have $\rho(f)\ge 0$ and $\rho(f)=0$ if and only if $L_\p/K_p$ is
%tame for all $\p|p$, that is if $p\nmid e_\p$ and the 
%residue field extension is separable.

\begin{proposition}\label{prop:deltastar}
  We have $\delta^*(f)\le \delta(f)-\rho(f)$.
\end{proposition}
\begin{proof} If $F\in \Oc_p[x]$ is monic irreducible with Okutsu
  frame $[\phi_1,\ldots,\phi_{r+1}]$, it follows from Definition
  \ref{def:okutsu} that
  $\delta_0(F)=\deg(F)v_p(\phi_{r}(\alpha))/\deg(\phi_{r})$ and
  \cite[Lemma 1.5]{BaNaSt13} shows that $\delta_0(F)$ coincides with
  the quantity introduced in \cite[Definition 2.1]{BaNaSt13}. We thus
  get $\delta_0(F)= \mu_r+\nu_r \le 2\mu_r=:2\mu(F)$ where $\mu_r$ and
  $\nu_r$ are Okutsu invariants defined by the successive slopes of
  the generalised Newton polygon encountered during the OM algorithm
  called with parameter $F$ (see \cite[page 141]{BaNaSt13}).
  % In particular, we have
  % $\mu_r=\mu(F):=\max \left\{v_p(g(\alpha)),\; g\in
  %   \Oc[x]\,\,\mathrm{monic},\, \deg(g)<\deg(F)\right\}$
  % \cite[Proposition 1.4]{Na14}.
  Combined with \cite[Proposition 1.4]{Na14}, we get for each $\p|p$
  \begin{equation}\label{eq:ddelta0}
    \frac{d_\p \delta_0(F_\p)}{2} \le d_\p \mu(F_\p) = \delta(F_\p)-f_\p \rho_\p
  \end{equation}
  where $\rho_\p\in \Ni$ is related to the local different by
  $ v_\p(\mathrm {Diff}(L_\p/K_p))=e_\p-1+\rho_\p.  $ Applying the
  norm $N_{L_\p/K_p}$ and using that
  $ \Delta_{L_\p/K_p}=N_{L_\p/K_p}(\mathrm {Diff}(L_\p/K_p)), $ we get
  $\sum_\p f_\p \rho_\p=\rho(f)$. The claim follows from summing
  \eqref{eq:ddelta0} over all $\p|p$ together with the classical
  formula
  $ \delta(f)=\sum_\p \delta(F_\p)+\sum_{\p\ne \q}
  v_p(\Res(F_\p,F_\q)).$ $\hfill\qed$
\end{proof}

The Okutsu invariant $\delta^*(f)$ is also closely related to the
$p$-index. In what follows, we use notations
\[
  \ind_p(f):=v_p([B:A[\theta]])\qquad \mathrm{and} \qquad
  \ind_p(F_\p)=v_p([B_\p:A_p[\theta_\p]]).
\]
Thus $\ind_p(f)=v_p(D)$ with notations of the introduction.

\begin{proposition}\label{prop:delta_index}
  We have the inequalities
  $ \mathrm{ind}_p(f) \le \delta^*(f)\le 2\,\mathrm{ind}_p(f)+d-1.$
\end{proposition}

\begin{proof}
  The first inequality follows from
  \[
    \frac{d_\p\delta_0(F_\p)}{2}\ge \frac{d_\p}{2}
    \mu(F_\p)=\ind_p(F_\p)+1-e_\p^{-1}\ge \ind_p(F_\p)
  \]
  (the equality by \cite[Proposition 1.4]{Na14}) together with
  \[
    \ind_p(f)=\sum_\p \ind_p(F_\p)+\frac{1}{2} \sum_{\p\ne \q}
    v_p(Res(F_\p,F_\q))
  \]
  (see e.g. \cite[Section 2.2]{Na14}). For the second inequality,
  \eqref{eq:disc} leads to
  \[
    \delta(f)=2\,\ind_p(f)+v_p(\disc[](L/K))=2\,\ind_p(f)+\sum_{\p|p}
    v_p(\disc[](L_\p/K_\p)).
  \]
  Combined with Proposition \ref{prop:deltastar}, we get
  $ \delta^*(f)\le 2\,\ind_p(f)+\sum_\p f_\p(e_\p-1) $ and we
  conclude thanks to the fundamental equality
  $\sum_{\p|p} e_\p f_\p=d_\p$. $\hfill\qed$
\end{proof}

\begin{theorem}\label{thm:precision}
  Let $\sigma\in \Ni$ and let $g,f\in \Oc[x]$ be two monic separable
  polynomials of degree $d$ such that $g\equiv f\mod p^{\sigma}$.
  \begin{enumerate}
  \item \label{enum:irr} If $\sigma > 2\delta^*(f)/d$, then $g$ is
    irreducible if and only if $f$ is irreducible.
  \item \label{enum:facto} If $\sigma > \delta^*(f)$, then any
    OM-factorisation of $g$ is an OM-factorisation of $f$.
  \end{enumerate}
\end{theorem}

\begin{proof}
  Point 1 follows from a closer look at the proof of
  \cite[Lemma 2.8]{BaNaSt13}. Namely, the quantities
  $u_{i,s}/e_0\cdots e_{i-1}$ that appear in the proof are upper
  bounded by $2\delta(F_s)/n_s$ (with $n_s=\deg(F_s)$) although
  \cite[Lemma 2.2]{BaNaSt13} allows to use the Okutsu bound to get a
  sharper inequality
  \[
    \frac{u_{i,s}}{e_0\cdots e_{i-1}} \le \delta_0(F_s) \le
    \frac{2\delta(F_s)}{n_s}.
  \]
  Hence, we may replace $\delta(F_s)$ by $n_s\delta_0(F_s)/2$ in
  \cite[inequalities (2.5) and (2.6)]{BaNaSt13}. By definition of
  $\delta^*(f)$, we get that $\delta(f)$ can be replaced by
  $\delta^*(f)$ in the upper bound of \cite[Lemma 2.8]{BaNaSt13}. We
  may thus also use $\delta^*(f)$ instead of $\delta(f)$ in
  \cite[inequality (2.8)]{BaNaSt13} of the proof of \cite[Theorem
  2.3]{BaNaSt13}, leading to Point \ref{enum:irr}. Similarly, we may
  replace $\delta(f)$ by $\delta^*(f)$ in \cite[Lemma 3.12]{BaNaSt13},
  from which it follows that we may replace $\delta(f)$ by
  $\delta^*(f)$ in \cite[Theorem 3.13]{BaNaSt13}, proving Point
  \ref{enum:facto}. $\hfill\qed$
\end{proof}

% This theorem shows that we can conduct all computations with
% precision $\delta^*(f)$ (\emph{i.e.} modulo $p^{\delta^*(f)+1}$) in
% order to compute an OM factorisation of $f$.

\begin{example}
  Consider $f=x^2+t^N x +t\in F_2[[t]][x].$
  We check that $\delta^*(f)=1$ while $\delta(f)=2N$ can take
  arbitrarily large values for the fixed degree $d=2$. By Theorem \ref{thm:precision}, the
  OM-factorisation of $f$ only depends on $f \mod p^2$, and we
  definitely don't want to work at precision $\delta(f)=2N$ for such a
  polynomial.
\end{example}

In terms of the index, Theorem \ref{thm:precision} together with
Proposition \ref{prop:delta_index} shows that we can work at precision
$2\,\ind_p(f)+d$ to get an OM factorisation of $f$. In fact, we can do
slightly better.
% Let $ \tilde{\delta}(f):=2\,\ind_p(f)+1=2\,v_p(D)+1.  $ We still
% denote $v_p$ the Gauss valuation on $K[x]$ attached to $v_p$.

\begin{theorem}\label{thm:precision_tilde_delta}
  Let $ \tilde{\delta}(f):=2\,\ind_p(f)+1.  $ Let
  $g\in \Oc[x]$ be a monic separable polynomial of degree $d$ such
  that $g\equiv f\mod p^{\sigma}$ for some
  $\sigma > \tilde{\delta}(f)$.
  \begin{enumerate}
  \item \label{enum:facto-prec} Any OM-factorisation of $g$ is an
    OM-factorisation of $f$.
  \item \label{enum:OM-prec} Running the OM algorithm of
    \cite{PoWe22b} with precision $\sigma$ returns an OM-factorisation
    of $f$ where the approximant $\phi_\p$ of $F_\p$ satisfy
    $ v_p(F_\p-\phi_\p)>\sigma-\ind_p(f).  $
  \item \label{enum:approx-prec} The approximants $\phi_\p$ satisfy
    the conditions of Corollary \ref{cor:triangular_p_basis} of
    Section \ref{sec:p_int_basis}.
  \end{enumerate}
\end{theorem}

\begin{proof} 
  The proof is algorithmic. We first call algorithm
  \texttt{Irreducible} of \cite{PoWe22b} with precision $\sigma$. By
  Proposition \ref{prop:delta_index}, we have
  $\tilde{\delta}(f)\ge 2\delta^*(f)/d$, so Theorem
  \ref{thm:precision} ensures that either we can certify $f$ that is
  irreducible (and compute in such a case an Okustsu approximant of
  $f$), or we detect a first partial factorisation that can be
  computed up to precision $\sigma$ thanks to a valuated Hensel lemma,
  getting
  \[
    f\equiv G_0\cdots G_r \mod p^{\sigma}.
  \]
  Each $G_i$ can be lifted to a factor $F_i$ of $f$ (not necessarily
  irreducible). Denoting $\mu$ the current augmented valuation
  (denoted $w$ in \cite{PoWe22b}) and $e$ the current ramification
  index (which is $w(\pi)$ in \cite{PoWe22b}) we deduce from
  \cite[Lemma 9]{PoWe22b} that
  \[
    \sigma_i:=v_p(G_i-F_i) \ge \sigma-\frac{\mu(G_i)}{e}.
  \]
  Now, denoting $\hat{F}_i$ the cofactor of $F_i$ in $f$, we deduce from \cite[Proposition 3.5]{BaNaSt13}
  $\frac{\mu(G_i)}{e}=\frac{\mu(F_i)}{e}=\frac{v(\Res(F_i,\hat{F}_i))}{deg(\hat{F}_i)}
  \le v(\Res(F_i,\hat{F}_i))$.
  \cite[Section 2.2]{Na14} leads to
  $\ind_p(f)=\ind_p(F_i)+\ind_p(\hat{F}_i)+ v(\Res(F_i,\hat{F}_i))$.
  As by assumption $\sigma > 2\,\ind_p(f)+1$, we get
  \[
    \sigma_i > \ind_p(f)+\ind_p(F_i)+\ind_p(\hat{F}_i)+1\ge
    \ind_p(f)+\ind_p(F_i)+1.
  \]
  As $\sigma_i> 2\,\ind_p(F_i)+1$, we can apply recursively this
  strategy on each approximant $G_i\equiv F_i\mod p^{\sigma_i}$,
  working now with precision $\sigma_i$. At a recursive call on an
  approximant $G$ of a factor $F$ of $f$, the current precision
  $\sigma'$ satisfies
  \[
    \sigma' \ge \sigma-\ind_p(f) > \ind_p(f)+\ind_p(F)+1\ge
    2\ind_p(F)+1\ge 2\delta^*(F)/\deg(F)
  \]
  so that the algorithm terminates and provides a complete
  OM-factorisation
  \[
    f\equiv \prod_\p G_\p \mod \pi^\sigma \quad \mathrm{with} \quad
    v_p(G_\p-F_\p)\ge \sigma-\ind_p(f)\quad \forall\,\, \p|p.
  \]
  This proves Points \ref{enum:facto-prec} and \ref{enum:OM-prec}. In particular, we get
  $v_p(G_\p-F_\p)> \ind_p(f)+1$. Using notations of Corollary
  \ref{cor:triangular_p_basis}, we have
  $\ind_p(f)=\lfloor \alpha_{1} \rfloor +\cdots \lfloor \alpha_{d-1}
  \rfloor \ge \lfloor \alpha_{d-1} \rfloor$ so that
  $v_p(G_\p-F_\p)> \alpha_{d-1}$, proving Point
  \ref{enum:approx-prec}. $\hfill\qed$
\end{proof}

\begin{example}
  For the example $f=x^q+t^{q} x +t\in \Fi_q[t][x]$ of the
  introduction of degree $d=q$ (a prime), we get $ \delta(f)=d^2$,
  $\delta^*(f)=d-1$ and $\tilde{\delta}(f)=1$ : we gain an extra
  factor $d$ for the precision when considering $\tilde{\delta}(f)$
  rather than $\delta^*(f)$.
\end{example}

\begin{remark}
  The bound $\sigma > 2\,\ind_p(f)+1$ is sharp at least when $\ind_p(f)=0$. Consider for instance
  $f=\prod_{i=0}^n ((x-i)^2-p)+p^N $ with $N>>0$, satisfying
  $\ind_p(f)=0$. Factoring $f\equiv \prod_i (x-i)^2\mod p $ would be
  neither sufficient to compute an OM-factorisation nor to detect if
  $\ind_p(f)=0$.  Factorisation modulo $p^2$ is required (see also Lemma \ref{lem:p_divides_D}).
  % Note as a remark that if we know by advance that $\ind_p(f)=0$,
  % then Dedekind's theorem ensures that the factorisation of $f$ mod
  % $p$ coincides with the reduction mod $p$ of the irreducible
  % factorisation of $f$.
\end{remark}

\begin{remark}
  A closer look at the proof shows that a precision
  $2\max \delta^*(g)/\deg(g)$ is sufficient for computing an
  OM-factorisation, where the max runs over all monic factors
  (possibly reducible) $g\in \Oc_p[x]$ of $f$. In most cases (when the index is not mainly due to one factor of small degree), this leads to a much smaller precision $\Oc(\delta^*(f)/d)$, which belongs to $\Oc(\ind_p(f)/d)$ by Proposition \ref{prop:delta_index}. This is the case for instance in the example detailled in Section \ref{sec:example}.
\end{remark}

%%% Local Variables:
%%% mode: latex
%%% TeX-master: "integral"
%%% End:

%% file: triangular.tex
We keep notations of Section \ref{sec:locbasis}.
% : $A$ is a PID with field of fraction $K$ and $L=K(\theta)$ is a
% field extension defined by the root of a separable monic irreducible
% polynomial $f\in A[x]$.
Let $B\subset L$ stands for the integral closure of $A$ in $L$. Denote
$A_p$ the localisation of $A$ at a fixed prime $p\in A$ and $B_p$ the
integral closure of $A_p$ in $L$. We have $B_p=B\otimes_A A_p$ and
$B_p$ is a free $A_p$-module of rank $d=\deg(f)$.

\begin{definition}
  A $p$-integral basis of $B/A$ (or $p$-basis) is an $A_p$-basis of
  $B_p$.
\end{definition}

There are several ways to compute a $p$-integral basis from the local
basis $\Bc_\p$ of the local rings $\Oc_\p$ for all $\p$ dividing $p$
\cite{BoDeLaPf22,GuMoNa13,GuMoNa15,Ba16}. Traditional methods (based
on the work of Okutsu \cite{Ok82}) compute a $p$-basis of shape
$\Bc_p:=\bigcup_{\p|p} z_\p \Bc_\p$ for some well chosen multipliers
$z_\p\in B_p$ (see also \cite{GuMoNa15} for the method of the
quotients). Although the complexity of such methods fit in our aimed
complexity bound, the resulting $p$-basis is not triangular in general
and the Hermite type reduction needed before applying CRT (Proposition
\ref{prop:CRT}) does not fit in our aimed complexity bound. The
wonderful \maxmin{} algorithm of Stainsby \cite{St18} avoids this
problem by providing directly a triangular $p$-integral basis.

% We have seen in the previous section how to compute an integral
% basis of $B_\p$ with respect to an Okutsu frame of $F_\p$. This
% Okutsu frame depend only on $F_\p$ truncated with finite
% $p$-precision and can thus be chosen in $A$. The corresponding
% integral basis of $B_\p$ then lies in $L$.

% Denote $e_\p=e(w_\p/v_p)$ the ramification index,
% $f_\p=[\Fi_\p:\Fi_p]$ and $d_\p=\deg(F_\p)$. We have
% $d_\p=e_\p f_\p$ and $d=\sum_{\p|p} e_\p f_\p$.

\subsection{Reduced triangular $p$-integral bases.}
For any prime ideal $\p$ dividing $p$, we define a valuation
$w_\p:L^\times \to \Qi$ by
\[
  w_\p(g(\theta)):=\frac{v_\p(g(\theta))}{e_\p}
\]
where $e_\p=e(\p/p)$ is the ramification index and $v_\p$ is the
canonical discrete valuation attached to $\p$. Thus $w_\p$ extends
$v_p$ to $L$ and $ w_\p(g(\theta))=v_p(g(\theta_\p)) $ where
$\theta_\p\in \overline{K}_p$ is any root of $F_\p$. Let $w=w_p$ be
the quasi-valuation defined by
\[
  w:L\to \Qi\cup\{\infty\},\qquad w(b):=\min (w_{\p}(b),\,\,\p|p).
\]
Thus, an element $b\in L$ belongs to $B_p$ if and only if $w(b)\ge 0$.
\begin{definition}\label{def:reduce}
  We say that a subset $\Bc=\{b_0,\ldots,b_{k}\}\subset L$ is
  $w$-\emph{reduced} if
  $ w(\sum_i \lambda_i b_i)=\min_i w(\lambda_i b_i) $ for all
  $\lambda_0,\ldots \lambda_k \in K$.
\end{definition}

Computing $w$-reduced integral bases is relevant for several
applications in function fields, such as computing Riemann-Roch spaces
\cite[Section 5]{Ba14} (possibly using various quasi-valuations
$w$). Given $g_0,\ldots,g_{k}\in A[x]$, we denote by
\[
  \Bc(g_0,\ldots,g_{k}):=\left(\frac{g_{0}(\theta)}{p^{ \lfloor
        w(g_0(\theta)) \rfloor}},\ldots,
    \frac{g_{k}(\theta)}{p^{\lfloor w(g_k(\theta)) \rfloor}}\right).
\]

\begin{definition} 
  A $p$-integral basis of shape $\Bc(g_{0},\ldots, g_{d-1})$ with
  $g_i$ monic of degree $i$ is called a triangular $p$-basis.
\end{definition}

For all $i=0,\ldots,d-1$, we define $\alpha_i=\alpha_i(f,p)\in \Qi^+$
by
\begin{equation}\label{eq:alpha_i}
  \alpha_i=\max\left\{w(h(\theta)),\,\, h\in A[x]\,\, \mathrm{monic},\,\,\deg(h)=i\right\}.
\end{equation}
We say that $g\in A[x]$ monic of degree $i$ is \emph{$w$-maximal} if
$w(g(\theta))=\alpha_i$.

\begin{proposition}\label{prop:p-triang-basis}\cite[Theorem 1.4]{St18}
  Let $\Bc=\Bc(g_{0},\ldots, g_{d-1})$, with $g_i$ monic of degree
  $i$. Then
  \begin{itemize}[label=$\bullet$]
  \item $\Bc$ is a $p$-integral basis if and only if
    $\lfloor w(g_i(\theta)) \rfloor=\lfloor \alpha_i\rfloor$.
  \item $\Bc$ is a reduced $p$-integral basis if and only if
    $w(g_i(\theta))=\alpha_i$.
  \end{itemize}
\end{proposition}

\begin{corollary}\label{cor:hi_versus_gi}
  Let $\Bc=\Bc(g_{0},\ldots, g_{d-1})$ and
  $\Bc'=\Bc(h_{0},\ldots, h_{d-1})$, with $g_i,h_i\in A[x]$ monic of
  degree $i$. If $\Bc$ is a triangular (resp. reduced triangular)
  $p$-basis and $h_i\equiv g_i\mod{} p^{\lfloor \alpha_{i} \rfloor }$
  (resp. $h_i\equiv g_i\mod{} p^{\lceil \alpha_{i} \rceil}$), then $\Bc'$ is a triangular (resp. reduced triangular) $p$-basis.
\end{corollary}

\begin{proof}
  Let $g,h\in A[x]$ and let $\alpha=w(g(\theta))$. If
  $h\equiv g \mod p^{\lfloor \alpha \rfloor }$ then
  $\lfloor w(h(\theta)) \rfloor =\lfloor \alpha \rfloor$ while if
  $h\equiv g \mod p^{\lceil \alpha \rceil}$ then
  $w(h(\theta))=\alpha$. Proposition \ref{prop:p-triang-basis}
  concludes. $\hfill\qed$
\end{proof}

By Proposition \ref{prop:p-triang-basis}, computing a reduced
triangular $p$-basis amounts to compute for each degree
$i=0,\ldots, d-1$ a maximal monic polynomial $g_i\in A[x]$ of degree
$i$.  It is shown in \cite{St18} that such polynomials can be obtained
as the product of exactly one Okutsu numerator of the local basis
$\Bc_\p$, for each $\p|p$.

\subsection{The \maxmin{} algorithm.}
For each $\p|p$, denote $d_\p=\deg(F_\p)$ and consider an extended set
of Okutsu $\p$-numerators of $F_\p$ (Definition \ref{def:okutsu_num})
\[
  \Nc_\p:=\Nc(F_\p)=\{g_{\p,0},\ldots,g_{\p,d_\p-1},g_{\p,d_\p}=\phi_\p\}.
\]
For each multi-index $\jb=(j_\p)_{\p|p}$ with $0\le j_\p \le d_\p$, we
define $g_\jb:=\prod_{\p|p} g_{\p,j_\p}$, so that $g_\jb\in A[x]$ is
monic of degree $\deg(\jb):=\sum_{\p|p} j_\p\le d$.

\begin{theorem}\label{thm:MaxPol}\cite[Theorem 2.6]{St18}
  If the approximants $\phi_\p$ are computed with a sufficient
  precision, then there exists for each $i=0,\ldots,d-1$ a multi-index
  $\jb_i$ of degree $i$ such that $g_i:=g_{\jb_i}$ is maximal,
  i.e. $w(g_i(\theta))=\alpha_i$.
\end{theorem}

We want to look for such optimal multi-indices without taking care of
the precision of the approximant $\phi_\p$. To this aim, we rather
consider $\phi_\p$ as a symbol and use the map
\[
  w_\p(g_\jb)=\left\{
    \begin{array}{ll}
      w_\p(g_\jb(\theta)) & \textrm{if } \phi_\p \, \nmid \, g_\jb\\%
      \,\infty & \textrm{if } \phi_\p \, | \, g_\jb.
    \end{array}
  \right.
\]
Let accordingly $w(g_\jb):=\min \{w_\p(g_\jb),\,\, \p|p\}$. We have
$w(g_\jb)<\infty$ if $\deg(\jb)<d$.

\begin{definition}
  We say that the multi-index $\jb$ is maximal if
  $w(g_\jb)\ge w(g_\ib)$ for all multi-index $\ib$ with
  $\deg(\ib)=\deg(\jb)$.
\end{definition}

Denote $\Pc=\{\p_1,\ldots,\p_s\}$ the set of prime ideals of $B$
dividing $p$ in a given fixed order. Denote $(\eb_1,\ldots,\eb_s)$ the
canonical basis of $\Zi^s$.
% Up to choosing a good ordering on the set $\Pc$, the following
% algorithm will provide a set of maximal multi-indices.

\begin{algorithm}[ht]
  \nonl\TitleOfAlgo{\maxmin{}($\Nc_{\p_1},\ldots,\Nc_{\p_s}$)\label{algo:MaxMin}}
  \KwIn{A set $\Nc_{\p_1},\ldots,\Nc_{\p_s}$ of local Okutsu
    numerators of $f$.}%
  \KwOut{Some maximal multi-indices $\jb_0,\ldots,\jb_{d-1}$ of
    degrees $0,\ldots,d-1$.}%
  $\jb_0\gets(0,\ldots,0)$\;%
  \For{$k=0,\ldots d-1$}{%
    $j\gets \min\{1\le i \le s,\,
    w_{\p_i}(g_{\jb_k})=w(g_{\jb_k})\}$\;%
    $\jb_{k+1}\gets\jb_k+\eb_j$\;%
  }%
\end{algorithm}

\begin{theorem}\label{thm:MaxMin}\cite[Theorem 3.3]{St18}
  There exists an ordering on the set $\Pc=\{\p_1,\ldots,\p_s\}$ such
  that algorithm \maxmin{} returns a correct answer.
\end{theorem}

We say in such a case that $\Pc$ is well-ordered. Such an order can be
read for free on the tree of types induced by the OM-factorisation of
$f$, see \cite[Section 3.2]{St18} for details.  Note that despite of
the remarkably simplicity of algorithm \maxmin{}, the proofs of
Theorems \ref{thm:MaxPol} and \ref{thm:MaxMin} in \cite{St18} are
quite involved.

% , they require in particular to compare carefully the various
% valuations $w_\q(\phi_{\p,i_\p}(\theta))$ thanks to closed formulas
% in terms of the Okutsu invariants computed in course of the OM
% algorithm .

\begin{corollary}\label{cor:triangular_p_basis}
  Let $\jb_0,\ldots,\jb_{d-1}$ be a sequence of maximal multi-indices,
  as returned by algorithm \maxmin{} (assuming $\Pc$
  well-ordered). Denote $g_i:=g_{\jb_i}$. If the precision of the
  approximant satisfies $v_p(F_\p-\phi_\p)\ge \alpha_{d-1}$ for all $\p|p$,
  then $g_i$ is a degree $i$ maximal polynomial for all
  $i=0,\ldots,d-1$. We have $w(g_{i}(\theta))=\alpha_i$ and get the reduced triangular $p$-basis
%  \begin{equation}\label{eq:triang_red_basis}
  \[
    \Bc=\left(1,\frac{g_{1}(\theta)}{p^{\lfloor \alpha_1 \rfloor }},\ldots, \frac{g_{d-1}(\theta)}{p^{\lfloor \alpha_{d-1} \rfloor }}\right).
  \]
%  \end{equation}  
\end{corollary}

\begin{proof}
  We have by the very definition
  \[
    w(g_i)=\min \{w_\p(g_i(\theta)),\,\, \p|p,\phi_\p\nmid g_i\}\ge
    \min \{w_\p(g_i(\theta)),\,\, \p|p\}=w(g_i(\theta)).
  \]
  If strict inequality holds then necessarily
  $w(g_i(\theta))=w_\p(g_i(\theta))$ for some prime $\p$ such that
  $\phi_\p$ divides $g_i$, from which it follows that
  $w(g_i(\theta))\ge w_\p(\phi_\p(\theta))$.  On the other hand, we
  have $\alpha_i\ge w(g_i(\theta))$ by definition of $\alpha_i$. We
  get
  $\alpha_i \ge w(g_i(\theta))\ge w_\p(\phi_\p(\theta)) =
  v(\phi_\p(\theta_\p))=v((\phi_\p(\theta_\p)-F_\p(\theta_\p))\ge
  v_0(F_\p-\phi_\p)\ge \alpha_{d-1}$. Since
  $\alpha_{d-1}\ge \alpha_i$, this forces $w(g_i(\theta))=\alpha_i$
  and we conclude with Proposition \ref{prop:p-triang-basis}.  Suppose
  now that equality $w(g_i)=w(g_i(\theta))$ holds. By Theorem
  \ref{thm:MaxPol}, there is some $g=g_\ib$ of degree $i$ such that
  $w(g(\theta))=\alpha_i$ (up to compute the $\phi_\p$'s with high
  enough precision). As $\jb_i$ is maximal, we have
  $w(g)\le w(g_i)$ independently of the chosen precision. We get
  $w(g_i(\theta))=w(g_i)\ge w(g)\ge w(g(\theta))=\alpha_i \ge
  w(g_i(\theta))$, the last inequality by definition of
  $\alpha_i$. Again, this forces
  $w(g_i(\theta))=\alpha_i$. $\hfill\qed$
\end{proof}

%%% Local Variables:
%%% mode: latex
%%% TeX-master: "integral"
%%% End:

%% file: global.tex
For each prime $p\in A$, we saw how to compute a (reduced)
$p$-integral basis
\[
  \Bc_p=\left(1,\frac{g_{p,1}(\theta)}{p^{\eta_{p,1}}},\ldots,
    \frac{g_{p,d-1}(\theta)}{p^{\eta_{p,d-1}}}\right)
\]
with $g_{p,i}\in A[x]$ monic of degree $i$ and
$\eta_{p,i}=\lfloor w_p(g_{p,i}(\theta))\rfloor$. We can glue these
various $p$-bases to get a triangular integral basis thanks to the
following result, due to Okutsu \cite[Thm 1]{Ok82IV} (see also
\cite[Thm 1.17]{St14} or \cite[Lem 1.3.18]{Ba14}):

\begin{proposition}
  \label{prop:CRT}
  Suppose given a $p$-basis $\Bc_p$ as above for each prime $p\,|\,
  D_f$. For all $i=0,\ldots,d-1$, let $h_i\in A[x]$ monic of degree $i$
  such that $h_i\equiv g_{p,i} \mod p^{\eta_{p,i}+1}$ for all
  $p\,|\,D_f$.  Then, the following family is a triangular integral basis of $L/K$:
  \[
    \Bc=\left(1,\frac{h_1(\theta)}{\prod_p p^{\eta_{p,1}}},\ldots,
      \frac{h_{d-1}(\theta)}{\prod_{p} p^{\eta_{p,d-1}}}\right).
  \] 
\end{proposition}

%Let $\Delta_{red}\in A$ be the radical of $\Delta_{L/K}$. We get as a
%corollary :
%
%\begin{theorem}\label{thm:same_integral_basis}
%  Let $g\in A[x]$ monic, irreducible of degree $d$ defining a field
%  extension $L'=K(\theta')$, and $B'$ the integral closure of $A$ in
%  $L'$. If $g\equiv f \mod (D\Delta_{red})^2$, then $B'=\tau(B)$,
%  where $\tau:K[\theta]\to K[\theta']$ is the $K$-vector space
%  isomorphism defined by $\tau(\theta)=\theta'$. In particular,
%  $D(g)=D(f)$.
%\end{theorem}
%
%\begin{proof}
%  The given condition is equivalent to
%  $g\equiv f \mod p^{2 v_p(D)+2}$ for all $p|\Delta(f)$. We conclude via Theorem
%  \ref{thm:precision_tilde_delta} combined with Corollary
%  \ref{cor:triangular_p_basis} and Proposition \ref{prop:CRT}.
%  $\hfill\qed$
%\end{proof}

%%% Local Variables:
%%% mode: latex
%%% TeX-master: "integral"
%%% End:

%% file: fractional.tex
\subsection{Fractional ideals}
We keep notations and hypothesis of the previous section.  Recall that
any fractional ideal $I$ of $B$ is a free $A$-module of rank $d$.

\begin{definition}\label{def:triangular_basis_fractional_ideal}
  A triangular basis (with respect to $f$) of a fractional ideal $I$
  of $B$ is an $A$-basis of $I$ of shape
\begin{equation}\label{eq:triang_basis_fract_ideal}
  \Bc_I=\left(\frac{1}{a_0},\frac{g_1(\theta)}{a_1},\ldots,\frac{g_{d-1}(\theta)}{a_{d-1}}\right),
\end{equation}
  where $a_i\in K^\times$ satisfy
  $a_{d-1} A\subset \cdots \subset a_0 A$ and $g_i\in A[x]$ is monic
  of degree $i$.
\end{definition}
Any fractional ideal admits a triangular basis \cite[Theorem
1.16]{St14}. The fractional ideals $a_i A$ of $A$ depend on the choice
of $f$ used to represent the field $L/K$ but the first fractional ideal
$a_0A$ does not:
\begin{lemma}\label{lem:IcapK}
  We have $a_0^{-1} A= I\cap K$.
\end{lemma}

\begin{proof}
  Let $\Bc_I=(b_0,\ldots,b_{d-1})$ be a triangular basis of $I$ and
  let $\alpha\in I$. Hence $\alpha=\sum \alpha_i b_i$ for some
  uniquely determined $\alpha_i\in A$. We have $\alpha\in K$ if and
  only if $\alpha$ has degree zero as a polynomial in $\theta$. Since
  $\Bc_I$ is triangular, this is equivalent to
  $\alpha_1=\cdots=\alpha_{d-1}=0$. Hence
  $I\cap K=b_0 A=a_0^{-1} A$.$\hfill\qed$
\end{proof}

\begin{definition}
  With notations as above, we define the normalised ideal of $I$ as
  $I^*:=a_0 I$.  We say that $I$ is normalised if $I^*=I$.
\end{definition}

Notice that $1\in I^*$ so that $B\subset I^*$. If $\Bc$ is a
(triangular) integral basis of $I^*$, then obviously $a_0^{-1}\Bc$ is a
(triangular) integral basis of $I$, and we may focus on the
computation of an $A$-basis of a normalised ideal.

In all what follows, we will assume that we are given the unique
factorisation of $I$ in terms of the prime ideals of $B$, denoted by
%\[
  $I=\prod_{\p} \p^{n_\p}$, $n_\p\in \Zi$.
%\]
There are efficient algorithms to determine such a factorisation given
the OM-factorisations of $f$ above the involved primes $p\in A$
\cite{GuMoNa13}.

The normalised ideal $I^*$ is then easily deduced.  Let $\Pc_A$ be the
set of primes of $A$ ({\it i.e.} a set of generators of the prime
ideals of $A$). For $p\in \Pc_A$, we define
\begin{equation}\label{eq:mp}
  m_p=m_p(I):=\max\left\{\left\lceil\frac{n_\p}{e_\p}\right\rceil,
    \,\,\p|p\right\} \, \in \Zi.
\end{equation}

\begin{lemma}\label{lem:valIpstar}
  We have $v_p(a_0)=-m_p$ and
  $ I^*=\prod_{p\in \Pc_A} \prod_{\p|p} \p^{n_\p-e_\p m_p}.  $
\end{lemma}

\begin{proof}
  As $B\subset I^*$, we have $v_\p(I^*)\le 0$ for all $\p$. This
  implies that $v_\p(a_0)\le -n_\p$. If $\p|p$, then
  $v_\p(a_0)=e_\p v_p(a_0)$, hence $-v_p(a_0)\ge n_\p/e_\p$. Since
  $v_p(a_0)\in \Zi$ we deduce that $-v_p(a_0)\ge m_p$. If strict
  inequality holds for some $p$, we get
  \[
    -v_p(a_0)\ge m_p+1\quad \Longrightarrow \quad -v_p(a_0)e_\p\ge
    n_\p e_\p+e_\p\quad \forall\,\, \p|p,
  \]
  and $-v_\p(p a_0)\ge n_\p$ for all $\p$. This implies
  $(a_0 p)^{-1}\in I\cap K$ in contradiction with Lemma
  \ref{lem:IcapK}. Hence $v_p(a_0)=-m_p$ for all $p$ and
  $v_\p(I^*)=n_\p-e_\p m_p$ for all $\p|p$.$\hfill\qed$
\end{proof}

\paragraph{Normalised size of fractional ideals.}
Given two free $A$-submodules $I,I'\subset L$ of rank $d$ with
respective $A$-basis $\Bc=(b_1,\ldots,b_d)$ and
$\Bc'=(b'_1,\ldots,b'_d)$, the transition matrix $T\in K^{d\times d}$
from $\Bc$ to $\Bc'$ is defined by
$ (b'_1,\ldots,b'_d)T=(b_1,\ldots,b_d).  $ If we change the $A$-basis
of $I$ or $I'$, then $T$ is multiplied by a matrix in $GL_d(A)$. Hence
the following definition makes sense :

\begin{definition}
  For two free $A$-submodules $I,I'\subset K$ of rank $d$, the index
  $[I':I]$ is the fractional ideal of $A$ generated by the determinant
  of the transition matrix from an $A$-basis of $I$ to an $A$-basis of
  $I'$. %For a given prime $p\in A$, the $p$-index of $I$ is $\mathrm{ind}_p(I):=v_p([I:A[\theta]])$.
\end{definition}

The index is multiplicative : $[I:I'']=[I:I'][I':I''].$ Moreover, if
$I'\subset I$, then $[I:I']\subset A$ and $I'=I$ if and only if
$[I:I']=A$. 

If $I$ has triangular basis as in Definition \ref{def:triangular_basis_fractional_ideal}, then $[I:A[\theta]]=(a_0\cdots a_{d-1})$ and by transitivity of the index, we get
  \[v_p([I^*:A[\theta]])=v_p(a_0\cdots a_{d-1})-d v_p(a_0) \,\ge 0,
\]
 positivity since $a_i A \subset a_0 A$. 
%Finally, we have $[B:I]=N_{L/K}(I)$ so that
%\[
%\ind_p(I)=\ind_p(f)-v_p(N_{L/K}(I)).
% \]
% since $[I:A[\theta]]=[I:B][B:A[\theta]]$.

%\begin{definition}
%  For $p\in \Pc_A$, we define the \emph{normalised $p$-index} of a
%  fractional ideal $I$ as $ \eta_p(I):=v_p\,(\,[I^*:A[\theta]]\,).  $
%\end{definition}

%The $p$-index of a normalised fractional ideal is easily deduced from a triangular basis:

\begin{definition}
  Suppose either $A=\Zi$ and $h(p):=\log\,|p|$, or $A=k[t]$ and
   $h(p):=\deg(p)$. The \emph{normalised size} of $I$ is
  \[
    h(I):=\sum_{p\in \Pc_A} v_p([I^*:A[\theta]])\,h(p)\,\, \in \Ni.
  \]
\end{definition}

\begin{lemma}\label{lem:compare_index}
  Denote $D=D_f$ as in \eqref{eq:disc}. Let $I$ be a fractional ideal
  of $B$. Then $ h(I) \ge h(D) $ and equality holds if and only if
  $I=\alpha B$ for some $\alpha\in K\setminus\{0\}$.
\end{lemma}

\begin{proof}
  We have $ [I^*:A[\theta]]=[I^*:B][B:A[\theta]]$ and
  $v_p([B:A[\theta]])=v_p(D)$ so that $v_p([I^*:A[\theta]])=v_p([I^*:B]) +v_p(D)$. Since $B\subset I^*$, we
  have $v_p([I^*:B])\ge 0$, leading to $h(I) \ge h(D)$. Equality is equivalent to that
  $v_p([I^*:B])= 0$ for all prime $p$, that is $[I^*:B]=A$. Since $B\subset I^*$, this is equivalent to that $I^*=B$, proving the last claim. $\hfill\qed$
\end{proof}

\subsection{$p$-bases of fractional ideals}
Let us fix $p\in A$ a prime. Denote $I_p:=I\otimes_A A_p$ the
localisation of $I$ at $p$. Note that
% \begin{equation}\label{eq:Ip}
\[
  I_p=\{b\in L,\,\,\, v_\p(b)\ge n_\p\,\,\, \forall\,\,\p | p\}
\]
%\end{equation}
and $I_p$ is a fractional ideal of $B_p$. As such, it is a free
$A_p$-module of rank $d$. A $p$-basis of $I$ is by definition an
$A_p$-basis of $I_p$. To compute such a basis, we can follow exactly
the same strategy than for the case $I_p=B_p$, except that we consider
now the shifted valuations
\[
  w_{\p,I}:L\to \Qi\cup\{\infty\},\qquad \quad
  w_{\p,I}(b(\theta)):=w_{\p}(b(\theta))- \frac{n_\p}{e_\p}
\]
and accordingly the map $w_I=w_{p,I}$ defined by
\[w_I:L\to \Qi\cup\{\infty\}, \qquad \quad w_I(b(\theta))=\min
  (w_{\p,I}(b(\theta)),\,\,\p|p).\] Thus, an element $b\in L$ belongs
to $I_p$ if and only if $w_I(b)\ge 0$.
Given $g_0,\ldots,g_{k}\in A[x]$, we denote by
\[
  \Bc_I(g_0,\ldots,g_{k}):=\left(\frac{g_{0}(\theta)}{p^{ \lfloor
        w_I(g_0(\theta)) \rfloor}},\ldots,
    \frac{g_{k}(\theta)}{p^{\lfloor w_I(g_k(\theta)) \rfloor}}\right).
\]

\begin{definition} 
  A triangular $p$-basis of $I$ is a basis of shape
  $\Bc_I(g_{0},\ldots, g_{d-1})$ with $g_i$ monic of degree $i$.
\end{definition}

We let
%\begin{equation}\label{eq:alpha_Ii}
  $\alpha_{I,i}=\max\left\{w(h(\theta)),\,\, h\in A[x]\,\, \mathrm{monic},\,\,\deg(h)=i\right\}$
%\end{equation}
for $i=0,\ldots,d-1$ and say that $g\in A[x]$ monic of degree $i$ is $w_I$-maximal if
$w_I(g(\theta))=\alpha_{I,i}$.

\begin{theorem}\label{thm:triangular_p_basis_fractional_ideal}
  Suppose $I$ normalised. Let $\jb_0,\ldots,\jb_{d-1}$ be a sequence
  of maximal multi-indices, as returned by Algorithm \maxmin{}
  called with the map $w_I$ instead of $w$ (assuming $\Pc$
  well-ordered). Let $g_i:=g_{\jb_i}$. If %the approximants satisfy
    $v_p(F_\p-\phi_\p)\ge \alpha_{I,d-1}$ for all $\p|p$,
  then $\Bc_I(g_0,\ldots,g_{d-1})$ is a triangular $p$-basis of $I$.
\end{theorem}

\begin{proof}
  It follows from respectively \cite[Theorems 1.25]{St14}, \cite[Theorem
  5.1]{St14} and \cite[Proposition 5.2]{St14} that the analogous of
  respectively Proposition \ref{prop:p-triang-basis}, Theorem
  \ref{thm:MaxPol} and Theorem \ref{thm:MaxMin} still hold if we
  replace $w$ by $w_I$ and $\alpha_{i}$ by $\alpha_{I,i}$. This mainly
  follows from the fact that the valuations $w_{\p,I}$ are simply a
  shift of the valuations $w_\p$. Since $I=\prod \p^{n_\p}$ is assumed
  to be normalised, we have $B\subset I$ hence $n_\p\le 0$ for all
  $\p$. Thus, \cite[Theorem 5.3]{St14} ensures that the analogous of
  Corollary \ref{cor:hi_versus_gi} holds too if we replace $w$ by
  $w_I$ and $\alpha_{i}$ by $\alpha_{I,i}$. The proof of Theorem
  \ref{thm:triangular_p_basis_fractional_ideal} is then \emph{mutatis
    mutandi} identical to the proof Corollary
  \ref{cor:triangular_p_basis}. $\hfill\qed$
\end{proof}

\subsection{Improvements \textit{via} $S$-basis}
Let $S\subset \Pc$. For $b\in L$, define $w_{I,S}(b)=\min
  (w_{\p,I}(b),\,\,\p\in S)$.
If we apply algorithm \maxmin{} with the set of denominators  $\{\Nc_\p,\p\in S\}$ as input and with $w_{I,S}$ instead of $w$, we get a family of multi-indices $\jb_0,\ldots,\jb_{d_S-1}$ with $\jb_i=(\jb_{i,\p})_{\p\in S}$, and where $d_S:=\sum_{\p\in S} \deg(F_\p)$.
The resulting polynomials $g_{i}:=g_{\jb_i}$ have maximal $w_{S,I}$-valuation (assuming that the involved $\phi_\p$ are computed with a sufficient precision) and give rise to a triangular set 
\[
  \Bc_{I,S}:=\left(\frac{g_{0}(\theta)}{p^{ \lfloor
        w_{I,S}(g_0(\theta)) \rfloor}},\ldots,
    \frac{g_{d_S-1}(\theta)}{p^{\lfloor w_{I,S}(g_{d_S-1}(\theta)) \rfloor}}\right)
\]
that we call \emph{an $S$-basis of $I$}. 
%which is a (reduced) triangular basis of the free $\Oc_p$-module $I_S:=\oplus_{\p\in S} \,\p^{n_\p}\,\Oc_\p$.  
Besides playing a key role in the proof of the Theorem \ref{thm:MaxMin} in \cite{St18}, these $S$-bases are also relevant to accelerate the computation of a triangular $p$-basis of $I$ in some particular cases. In what follows we let 
\begin{equation}\label{eq:T}
T=\left\{\p\in \Pc\,,\,\,  \mathrm{ind}_p(F_\p)=v_p(\mathrm{Res}(F_\p,\hat{F}_\p))=n_\p=0\right\}
\end{equation}
and we denote $S=\Pc\setminus T$.

% , and $w_\p(\phi_\q)=0$ for all $\p\ne \q$ such that $F_\p|g$ and
% all $F_\q|h$ or $F_\p|h$ and $F_\q|g$, or $F\p|g$ and
% $F_\q|g$. Also, $w_\p(\phi_\p)=0$ for all $F_\p|g$, except possibly
% if $g$ has a factor of shape $x-\alpha$ with $v_p(\alpha)>0$.

\begin{proposition}\label{prop:first_facto}
   Suppose $I$ normalised. Let $g=\prod_{\p\in T}\phi_\p \in A[x]$ and consider $\Bc_{I,S}=(b_0,b_1,\ldots,b_{d_S-1})$ an $S$-basis of $I$ as above.  The set
 \[
    \Bc_I=(1,\theta,\ldots,\theta^{d-d_S-1},g\,b_0,,\ldots, g\, b_{d_S-1})).
\]
is a triangular $p$-integral basis of $I$.
\end{proposition}

\begin{proof} Note first that $\deg(g)=d-d_S$ so the set $\Bc_I$ is
  indeed triangular. By Proposition \ref{prop:p-triang-basis} (in the context of fractionary ideals), we need to show that the polynomials $x^i$ and $g b_j$  are $\lfloor w_I \rfloor$-maximal for $0\le i<\deg(g)$ and $0\le j<d_S$. Let $b$ be a $w_I$-maximal polynomial of degree $k<d$. By Theorem \ref{thm:MaxPol} (which remains valid with $w_I$ instead of $w$), we may take $b=\prod_{\p\in \Pc} b_\p$, with $b_\p\in \Nc_\p$. 
  
$\bullet$ If there exists $\p\in T$ such that $b_\p\ne \phi_\p$, then  $\deg(b_\p)<d_\p$, which forces $\lfloor w_\p(b_\p) \rfloor\le \ind_p(F_\p)$ (Proposition \ref{prop:local_basis}), hence $\lfloor w_{I,\p}(b_\p)\rfloor=0$ since $\ind_p(F_\p)=n_\p=0$. But $v_p(\mathrm{Res}(F_\p,\hat{F}_\p))=0$ implies also that $w_{I,\p}(b_\q)=0$ for all $\q\ne \p$. Hence $\lfloor w_{I,\p}(b)\rfloor=0$. As $I$ is normalised, we have $\lfloor w_{I,\q}(b)\rfloor\ge 0$ for all $\q$ and we deduce $\lfloor w_I(b)\rfloor=0$. As $b$ is $w_I$-maximal, we deduce that $\lfloor \alpha_{I,k}\rfloor=0$ and any monic degree $k$ polynomial (in particular $x^k$) is $\lfloor w_I \rfloor$-maximal.

$\bullet$ If $b_\p=\phi_\p$ for all $\p\in T$, then $b=g b'$ with $b'=\prod_{\p\in S} b_{\p}$. In particular, we have $k\ge \deg(g)$. We get \[w_{I}(b)=\min_{\p\in \Pc} (w_{I,\p}(b')+w_\p(g))=\min_{\p\in S} (w_{I,\p}(b')+w_\p(g))=\min_{\p\in S} w_{I,\p}(b')=w_{I,S}(b'),\]
the second equality because $w_\p(g)=\infty$ for $\p\in T$ and the third equality because $v_p(\mathrm{Res}(F_\q,\hat{F}_\q))=0$ for all $\q\in T$ forces $w_{\p}(g)=0$ for all $\p\in S$. Hence $b$ is $w_I$-maximal if and only if $b'$ is $w_{I,S}$-maximal and the claim follows.  $\hfill\qed$
\end{proof}

%\begin{proof}
%  Since $\ind_p(g)=0$, the irreducible factors $\phi_\p$ of $g$ are
%  pairwise coprime mod $p$, and since
%  $v_p(\Res(g,h))=0$, they are also coprime to $h$ mod $p$. This
%  implies that for $\alpha_i$ defined in \eqref{eq:alpha_i}, we get
%  $\alpha_i(f)=0$ for $i=0,\ldots,\deg(g)-1$, and
%  $\alpha_i(f)=\alpha_i(h)$ for $i=\deg(g),\ldots,d-1$. Finally, we
%  check that each numerator of $\Bc$ obeys to the conditions of
%  first item in Proposition \ref{prop:p-triang-basis} (see also
%  \cite[Proposition 6.1]{BoDeLaPf22}). $\hfill\qed$
%\end{proof}

\begin{remark} The basis $\Bc_I$ of Proposition \ref{prop:first_facto} is not necessarily $w_I$-reduced (Definition \ref{def:reduce}). Consider for instance $f=(x-1)^2+p\in A[x]$ and $I=B$. We have $\ind_p(f)=0$ and using Proposition \ref{prop:first_facto} would lead to the $p$-integral basis $\Bc=(1,\theta)$. This basis is not $w$-reduced as $w(\theta)=0< w(\theta-1)=1/2$. Using Corollary \ref{cor:triangular_p_basis} would have returned the \emph{reduced} $p$-basis $(1,\theta-1)$.
\end{remark}

\subsection{Global triangular bases of fractional ideals}

\begin{proposition}\label{prop:CRT_fractional}
  Let $I$ be a normalised fractional ideal of $L$. Suppose given a
  triangular $p$-basis $\Bc_{I,p}=\Bc_{I,p}(g_{p,0},\ldots,g_{p,d-1}) $
  of $I$ for each prime $p$ dividing $[I:A[\theta]]$ and
  let $\eta_{p,i}(I):=\lfloor w_{I}(g_{p,i}(\theta))\rfloor$. For
  all $i=0,\ldots,d-1$, let $h_i\in A[x]$ monic of degree $i$
  s.t. $h_i\equiv g_{p,i} \mod p^{\eta_{p,i}(I)+1}$ for all $p$
  dividing $[I:A[\theta]]$. Then a triangular $A$-basis of $I$ is given by
  \[
    \Bc_I=\left(1,\frac{h_1(\theta)}{\prod_p
        p^{\eta_{p,1}(I)}},\ldots, \frac{h_{d-1}(\theta)}{\prod_{p}
        p^{\eta_{p,d-1}(I)}}\right).
  \]
\end{proposition}

\begin{proof}
  This follows from \cite[Theorem 1.27]{St14}, where we use moreover that
  $g_{p,0}=1$ for all $p$ since $I$ is assumed to be
  normalised. $\hfill\qed$
\end{proof}

\begin{remark}
  If $I$ is not normalised, we first compute $I^*$
  following Lemma \ref{lem:valIpstar} and then compute an integral basis
  $\Bc_{I^*}$ of $I^*$ following Proposition
  \ref{prop:CRT_fractional}. An integral basis of $I$ is then given by
  $\Bc_I=\alpha \Bc_{I^*}$, with
  $\alpha:=\prod_{p|[I^*:A[\theta]]} p^{m_p}$, the integer $m_p$ being
  defined in \eqref{eq:mp}.
\end{remark}

\begin{theorem}\label{thm:main_fractional}
  Let $I$ be a fractional ideal. Given a squarefree factorisation of
the ideal  $[I^*:A[\theta]]$, we can compute a triangular integral basis of $I$
  with
  \begin{enumerate}
  \item $\Oc_\eps(\dx \,h(I))$ operations in $k$ if $A=k[t]$ with
    $\char(k)=0$ or $\char(k)>\dx$,
  \item $\Oc_\eps(\dx \,h(I)+h_D^2)$ operations in $k$ if $A=k[t]$
    with $\char(k)\leq\dx$,
  \item $\Oc_\eps(\dx\,h(I)+h_D^2)$ word operations if $A=\Zi$.
  \end{enumerate}
\end{theorem}
This result has to be compared to
$\Oc_\eps(d^3\,h(I)^2+d\,h_\Delta^2)$ that can be deduced from \cite[Thm
5.3.19]{Ba14}. The proof will be given at the end of Section \ref{sec:complexity}.

\begin{remark} Note that when $I=B$, we get $h(B)=h_D$ so the first statement of Theorem \ref{thm:main_hD} is a particular instance of Theorem \ref{thm:main_fractional}.
\end{remark}

\begin{remark}
  We don't take into account the last multiplications inherent to
  the relation $\Bc_I=\alpha\Bc_{I^*}$ in our estimation cost. 
  For many purposes, this step would lose crucial information.
\end{remark}

%%% Local Variables:
%%% mode: latex
%%% TeX-master: "integral"
%%% End:

%% file: complexity.tex
We keep notations and hypothesis of previous sections. In particular, $f\in A[x]$ is an irreducible separable monic degree $d$ polynomial.
For a prime $p\in A$, we denote $k_p$ the residue field of $\Oc_p$ and we charge one operation in $k_p$ for one operation in a fixed set $\Ac\subset A$ of representatives of $k_p$.

 \paragraph{Cost of the OM-factorisation.} 

%In order to use the complexity results of
%\cite{PoWe22b} with such a bound for the precision, we need to express
%the upper bound of \cite[Lemma 4]{PoWe22b} in terms of the $p$-index,
%i.e. to bound the residual degree $f_\p$ of each irreducible factor
%$F_\p$ of $f$.
%\begin{lemma}\label{lem:bound_for_f}
%  Let $F_\p\in \Oc_p[x]$ be monic irreducible, with residual degree
%  $f_\p$. If $\mathrm{ind}(F_\p)>0$, then $f_\p\le 2\,\mathrm{ind}(F_\p)$.
%\end{lemma}
%\begin{proof}
%  From \cite[Propoposition 1.4]{Na14}, we get
%  $d_\p\,\mu(F_\p) =2\,\ind(F_\p)+d_\p-f_\p$ (using $d_\p=e_\p\,f_p$),
%  i.e. $f_\p = 2\ind(F_\p) + d_\p(1-\mu(F_\p))$. We are done if
%  $\mu(F_\p)\ge 1$. From the formula for $\mu(F)$ in
%  \cite[Propoposition 1.4]{Na14} again, we see that $\mu(F_\p)=0$ only
%  if $e_\p=1$, thus $\ind(F_\p)=0$.
%\end{proof}
%In particular, for each $\p|p$ with $\ind(F_\p)>0$, we get
%$f_\p\le 2\,\ind_p(f)$. This rough estimate is enough for our
%purpose. Finally, when $\ind(F_\p)=0$ (that we can detect on
%$f \mod p^2$), there is no need to use \cite[Lemma 4]{PoWe22b} in the
%complexity analysis. 

\begin{proposition}\label{prop:comp-OM}
  Assume $\sigma\geq 2\,\mathrm{ind}_p(f)+1$. We
  can compute an OM factorization of $f$ above $p$ such that $v_p(\phi_\p-F_\p)\ge \sigma-\mathrm{ind}_p(f)$  for all $\p|p$ with $\Oe(\dx\,\sigma)$ operations in $k_p$ if $\char(k_p)=0$ or $>d$, and  $\Oe(\dx\,\sigma+\mathrm{ind}_p(f)^2)$ operations in $k_p$ otherwise. 
  The algorithm returns as a byproduct the values
  $w_\q(\phi_{\p,i}(\theta))$ for all $\p,\q|p$ and all
  $0\le i\le r_\p+1$.
\end{proposition}

\begin{proof}
When $\char(k_p)$ is zero or $>d$, this follows from Theorem
\ref{thm:precision_tilde_delta} together with \cite[Thm 4]{PoWe22b}. When $0<\char(k_p)<d$, the number of refinement steps inherent to the OM algorithm is bounded by $\Oc(\ind(f)/\ell_0)$ thanks to \cite[Def 4.15 and Thm 4.18]{GuMoNa12}, with $\ell_0$ being the first residual degree as used in \cite[Lemma 4]{PoWe22b}. We thus need to add a cost $\Ot(\delta^* \, \ind_p(f))$  thanks to \cite[Section 3]{PoWe22b} (replacing again $\delta$ by $\delta^*$ in \cite[Thm 2]{PoWe22b}). By Proposition \ref{prop:delta_index}, this is $\Ot(\ind_p(f)^2+d\,\ind_p(f))$ which fits in the aimed bound. $\hfill\qed$
\end{proof}

Note that we only compute the square-free factorisation of the various residual polynomials and we rely on dynamic evaluation, using the
complexity results of \cite{HoLe19} in that context (see \cite[Section 5.4]{PoWe22b}).

\paragraph{Binary cost of the \maxmin{} algorithm.}  We consider the general context of a normalised fractional ideal $I$ of $B$. 
Denote for short $\ind_p(I)=v_p(I:A[\theta])$. Note that $\ind_p(B)=\ind_p(f)=v_p(D)$. 

\begin{proposition}\label{prop:comp-maxmin}
  The cost of \maxmin{} above a prime $p$ with respect to the quasi-valuation $w_I$ is $\Oc(ds\log(\mathrm{ind}_p(I)))$ word
  operations, with $s$ the number of irreducible factors of $f$ in $\Oc_p[x]$.
\end{proposition}

\begin{proof}
  There are $d$ iterations including one minimum of a set of cardinality $s$ and one addition, each element having binary size  bounded by $\log(\mathrm{ind}_p(I))$. $\hfill\qed$
\end{proof}

With regards to Theorem \ref{thm:main_fractional}, we need to take care that we might have $s>\ind_p(I)$.  
 For such a small $p$-index, we rather use Proposition \ref{prop:first_facto}.

\begin{lemma}\label{lem:maxmin}
  Suppose $I$ normalised and let $S\subset \Pc$ as defined by \eqref{eq:T}. The binary cost of \maxmin{} to compute an $S$-basis of $I$ is $\Ot(d\,\mathrm{ind}_p(I))$.
\end{lemma}

\begin{proof}
  The cost of \maxmin{} restricted to $S$ is now $\Oc(d_S \,\Card(S)\, \log(\mathrm{ind}_p(I)))$. We have $\mathrm{ind}_p(I)=v_p([B:A[\theta])+v_p([I:B])=\ind_p(f)-v_p(N_{L/K}(I))$, leading to
  \[\ind_p(I)= \sum_{\p|p} \left(\ind_p(F_\p)+\frac 1 2 v_p(\Res(F_\p,\hat{F}_\p)-f_\p n_\p \right).
  \]
  Since $I$ is normalized, we have $n_\p\le 0$  for all $p$ so each summand is $\ge 1/2$ whenever $\p \in S$ by \eqref{eq:T}. We get $\Card(S)\le \ind_p(I)/2$ and the claim follows. $\hfill\qed$
\end{proof}

\paragraph{Cost of expanding and gluing $p$-integral basis.}

\begin{proposition}\label{prop:comp-local}
  Suppose $I$ normalised. Up to the cost of the OM-factorisation, one can
  compute a triangular $p$-basis of $I$ in less than $\Ot(d\,\mathrm{ind}_p(I))$ operations in $k_p$.
\end{proposition}

\begin{proof}
  Up to use Proposition \ref{prop:first_facto} and Lemma \ref{lem:maxmin}, we can compute non expanded denominators $g_0,\ldots,g_{d-1}$ of a triangular $p$-basis $\Bc_I$ in the aimed cost (binary cost for this step). There remains to expand $g_i$ mod $p^{\lfloor \alpha_{I,i}\rfloor}$ (Corollary \ref{cor:hi_versus_gi} and Corollary \ref{cor:triangular_p_basis} in the context of triangular ideal), for a cost of $\Ot(\deg(g_i)\lfloor \alpha_{I,i}\rfloor)$ operations in $k_p$. As $\sum_i \lfloor \alpha_{I,i}\rfloor = \ind_p(I)$ and $\deg(g_i)=i\le d$, the total cost is $\Ot(d\,\ind_p(I))$ operations in $k_p$.  $\hfill\qed$
\end{proof}

%\begin{corollary}\label{cor:intbasis-local}
%  We can compute an integral basis of $A_p[\theta]$ in a number of
%  operations in the residue field of $A_p$ bounded by
%  $\Oe(\dx\,\ind_p(f))$ if $\char(A_p)$ is zero or greater than $d$,
%  and $\Oe(\dx\,\ind_p(f)+\ind_p(f)^2)$ otherwise.
%\end{corollary}
%
%\begin{proof}
%  From Theorem \ref{thm:precision_tilde_delta}, one can run the
%  OM-algorithm with precision $\sigma \geq 2\,\ind_p(f)+1$. The result
%  is then a straight consequence of the propositions
%  above. 
%\end{proof}

\begin{proposition}\label{prop:comp-global}
Suppose $I$ normalised. Given a triangular $p$-integral basis above each prime $p\in A$ dividing $[I:A[\theta]]$, we can compute a global integral basis of $I$ in less than $\Ot(d h(I))$ binary operations if $A=\Zi$ or $\Ot(d h(I))$ operations in $k$ if $A=k[t]$. 
\end{proposition}

\begin{proof}
Let us first suppose that $A=k[t]$. Computing the polynomial $h_i\in A[x]$ in Proposition \ref{prop:CRT_fractional} requires $\Ot(\deg(h_i)\sum_p \eta_{p,i}(I) h(p))$ operations in $k$ by fast Chinese multi-remaindering. The result follows by summing over all $i=0,\ldots, d-1$, using $\deg(h_i)=i\le d$, $\sum_{i=0}^{d-1} \eta_{p,i}(I)=\ind_p(I)$ and $\sum_{p} \ind_p(I) h(p)=h(I)$. The same reasonning applies if $A=\Zi$, counting now the number of word operations. $\hfill\qed$
\end{proof}

\paragraph{Proof of Theorem \ref{thm:main_fractional}.} We may assume $I=I^*$. By Theorem \ref{thm:triangular_p_basis_fractional_ideal}, it's enough to compute the $\phi_\p$'s with precision $\mathrm{ind}_p(I)$. By Proposition \ref{prop:comp-OM}, we may apply the OM algorithm with precision $\sigma=\ind_p(I)+\ind_p(f)+1$, which costs $\Oe(d\,\ind_p(I))$ operations in $k_p$ (recall that $\mathrm{ind}_p(I)\ge \ind_p(f)$ since $I$ is normalized), plus an extra $\Oe(\ind_p(f)^2)$ for small residual characteristic.
The result then follows from Proposition \ref{prop:comp-OM}, Proposition \ref{prop:comp-local} and Proposition \ref{prop:comp-global}. Note that we rely again on dynamic evaluation since we are working above a squarefree factor $p$ of the index of $I$ which is not necessarily irreducible. $\hfill\square$

\paragraph{Proof of Theorem \ref{thm:main_hD}.} 
First part of Theorem \ref{thm:main_hD} follows from Theorem \ref{thm:main_fractional} applied with $I=B$. 
Let us prove the last claim, assuming now that we only know a squarefree factorisation of the discriminant $\Delta$. We use the following lemma.

\begin{lemma}\label{lem:p_divides_D}
Given a prime $p$, the condition $p|D$ only depends on $f\mod p^2$, and can be checked with $\Ot(d)$ operations in $k_p$. 
\end{lemma}

\begin{proof}
We want to check if $\ind_p(f)=0$. We first compute the square-free factorisation of $f\mod p$ and lift it once to get
$f=\prod_i F_i\mod p^2$, where $F_i|f$ and the $F_i$'s are coprime mod $p$. This costs $\Ot(d)$ operations in $k$. We have $\ind_p(f)=0$ if and only if $\ind_p(F_i)=0$ for all $i$. Since $F_i=P_i^{N_i}\mod p$ with $P_i\in A[x]$ monic and square-free mod $p$, we have $\ind_p(F_i)=0$ if and only if $N_i=1$ or $v_p(F_i \mod P_i)=1$ (Eisenstein case), see e.g. \cite[Rem 4.13]{GuMoNa12}. Both conditions only depend on $F_i\mod p^2$ and can be checked with $\Ot(\deg(F_i))$ operations in $k_p$ for each $i$, hence a total cost $\Ot(d)$.  $\hfill\qed$
\end{proof}

The last claim in Theorem \ref{thm:main_hD} follows. By \eqref{eq:disc}, we need to check if $p|D$ only if $p^2|\Delta$. By Lemma \ref{lem:p_divides_D}, this adds an extra cost of $\Ot(d\,\sum_{p^2|\Delta}h(p))= \Ot(d h_{red})$ (binary cost if $A=\Zi$ or arithmetic cost if $A=k[t]$). $\hfill\square$

%% file: stainsby-3.3.tex
We conclude our paper by illustrating the different steps of our
algorithm on the following example of Stainsby \cite[Section
3.3]{St18} over $\Zi[x]$, that is the polynomial
\[
  \begin{array}{@{}l@{}}
    f=x^{13}+{\scriptstyle 3}\,q^8\,x^{11}+{\scriptstyle 18753}\,q^{12}\,x^{10}+{\scriptstyle 781253}\,q^{16}\,x^9+{\scriptstyle 244178131}\,q^{20}\,x^8+\\%
    {\scriptstyle 783631254}\,q^{24}\,x^7+{\scriptstyle 14894940628}\,q^{28}\,x^6+{\scriptstyle 763967225003}\,q^{32}\,x^5+{\scriptstyle 193053764471876}\,q^{36}\,x^4\\%
    +{\scriptstyle 1562575008}\,q^{48}\,x^3+{\scriptstyle 488318756}\,q^{52}\,x^2+{\scriptstyle 1527929762506}\,q^{56}\,x+{\scriptstyle 4579209021877}\,q^{60}
  \end{array}
\]
with $q=5$. We have $\Delta_f=2^6\,5^{744}\,n$ with $n$ squarefree,
and $D_f=2^3\,5^{372}$. Knowing either of them, we first compute
integral bases over $p=2$ and $p=5$.

\paragraph{Reduced triangular basis over $p=5$.}
The OM algorithm run with precision $69$ finds a factorisation
$f=f_1\,f_2\,f_3$, for which we have initial approximations
\[
  \psi_1= x^4+2\,p^{24},
  \psi_2=\phi_1^2+p^{18}\,\phi_1+p^{32}\,x+p^{36},
  \psi_3=\phi_1+p^{17}
\]
with $\phi_1=x^3+p^8\,x+p^{12}$. It provides the three associated
quasi-valuations $w_1,w_2,w_3$ satisfying (we denote
$\war{}(a)=(w_1(a),w_2(a),w_3(a))$):
\[
  \begin{array}{c}
  \war{}(x) = (6,4,4)\ ;\ \war{}(\phi_1)=(12,18,17)\\%
  \war{}(f_1) = (\infty,16,16)\ ;\ \war{}(f_2) = (24,\infty,34)\ ;\
  \war{}(f_3) = (12,17,\infty)
  \end{array}
\]
We also get $w_1(\psi_1)=24$, $w_2(\psi_2)=36$, $w_3(\psi_3)=17$ and
the numerator sets are
\[
  \Nc_{5,1}=\{1,x,x^2,x^3,f_1\} ;
  \Nc_{5,2}=\{1,x,x^2,\phi_1,x\,\phi_1,x^2\,\phi_1,f_2\} ;
  \Nc_{5,3}=\{1,x,x^2,f_3\}.
\]
\maxmin{} runs as follows: starting from $g_0=1$, at each step $i$, we
increase in the triplet the index equal to the smallest
$j\in\{1,2,3\}$ s.t. $w_j(g_i)=\min w(g_i)$ (which is underlined in
the last column).
\[
  \begin{array}{c|c|c|c|c}
    i & \textrm{triplet} & g_i &  \war{}(g_i) &w(g_i)\\%
    \hline
    0 & (0,0,0) & 1   & (\underline 0,0,0) & 0\\%
    1 & (1,0,0) & x   & (6,\underline 4,4) & 4\\%
    2 & (1,1,0) & x^2 & (12,\underline 8,8) & 8\\%
    3 & (1,2,0) & x^3 & (18,\underline{12},12) & 12\\%
    4 & (1,3,0) & x\,\phi_1 & (\underline{18},22,21) & 18\\%
    5 & (2,3,0) & x^2\,\phi_1 & (\underline{24},26,25) & 24\\%
    6 & (3,3,0) & x^3\,\phi_1 & (30,30,\underline{29}) & 29\\%
    7 & (3,3,1) & x^4\,\phi_1 & (36,34,\underline{33}) & 33\\%
    8 & (3,3,2) & x^5\,\phi_1 & (42,38,\underline{37}) & 37\\%
    9 & (3,3,3) & x^3\,\phi_1\,f_3 & (\underline{42},47,\infty) & 42\\%
    10 & (4,3,3) & \phi_1\,f_1\,f_3 & (\infty,\underline{51},\infty) & 51\\%
    11 & (4,4,3) & x\,\phi_1\,f_1\,f_3 & (\infty,\underline{55},\infty) & 55\\%
    \;12\; & \;(4,5,3)\; & \;x^2\,\phi_1\,f_1\,f_3\; & \;(\infty,\underline{59},\infty)\;&\;59\;%
  \end{array}
\]
Note that any order on the prime ideals $\p_1$, $\p_2$, $\p_3$ works
here. We get:
% \[
%   \begin{array}{rl}
%     \Nc_5=\{&1,x,x^2,x^3,x\,\phi_1,x^2\,\phi_1,x^3\,\phi_1,x^4\,\phi_1,x^5\,\phi_1,x^3\,\phi_1\,f_3,
%               \phi_1\,f_1\,f_3,\\%
%             &x\,\phi_1\,f_1\,f_3,x^2\,\phi_1\,f_1\,f_3\}
%   \end{array}
% \]
\[
    \Nc_5=\{{\scriptstyle{} 1,x,x^2,x^3,x\,\phi_1,x^2\,\phi_1,x^3\,\phi_1,x^4\,\phi_1,x^5\,\phi_1,x^3\,\phi_1\,f_3,
              \phi_1\,f_1\,f_3,x\,\phi_1\,f_1\,f_3,x^2\,\phi_1\,f_1\,f_3}\}
\]
% l.9  3*w3(x)+w3(phi)+w3(psi3) = 29+w3(psi3) -> 13 ok
% l.10: w3(phi1)+w3(f1)+w3(psi3) = 33+w3(psi3) -> 18 ok ; w1(phi1)+w1(psi1)+w1(f3) = 24+w1(psi1) -> 28 ok (27 aussi ?)
% l.11: w3(x)+w3(phi1)+w3(f1)+w3(psi3) = 37+w3(psi3) -> 18 ok ; w1(x)+w1(phi1)+w1(psi1)+w1(f3) = 30+w1(psi1) -> 26 ok (25 aussi ?)
% l.12: 2*w3(x)+w3(phi1)+w3(f1)+w3(psi3) = 41+w3(psi3) -> 18 ok ; 2*w1(x)+w1(phi1)+w1(psi1)+w1(f3) = 36+w1(psi1) -> 24 ok (23 aussi ?)
Finally, one can check that it is sufficient to get approximations
$\psi_i$ of the factors $f_i$ satisfying $w_1(\psi_1)\geq 28$,
$w_2(\psi_2)\geq 36$ (that we got from the OM algorithm with precision
$69$) and $w_3(\psi_3) \geq 18$.  By lifting the factorisation only
once (using precision $70$), we get that, and update:
\[
  \psi_1=x^4+2\,p^{24}+p^{34} \;,\
  \psi_3=\phi_1+p^{17}+2\,p^{12}\,x^2+4\,p^{16}\,x+3\,p^{20}
\]

\paragraph{Reduced triangular basis over $p=2$.} Here we do not need
the \maxmin{} algorithm. As
$f\mod 2=(y^3+y^2+1)\,(y^4+y^3+1)\,(y^3+y+1)^2$, we only use the OM
algorithm above the factor $(y^3+y+1)$, then use Proposition
\ref{prop:first_facto} to conclude. The OM algorithm provides a local
set of numerators equal to $\{1,x,x^2,\phi_0,x\,\phi_0,x^2\,\phi_0\}$
with $\phi_0=x^3+x+1$. We thus get
\[
  \Nc_2= \{
  1,x,x^2,x^3,x^4,x^5,x^6,g,x\,g,x^2\,g,\phi_0\,g,x\,\phi_0\,g,x^2\,\phi_0\,g
  \}
\]
with $g=(x^3 + x^2 + 1) \, (x^4 + x^3 + 1)$.

\paragraph{Global basis.} We finally glue together these two bases
using CRT.  For instance, we compute $g_4$ s.t.
$g_4 = x^4+5^8\,x^2+5^{12}\,x \mod 5^{19}$ and $g_4 = x^4\mod 2$, that
is
$g_4=5^{19}\,x^4-2\frac{5^{19}-1}2(x^4+5^8\,x^2+5^{12}\,x)=x^4+(5^8-5^{27})\,x^2+(5^{12}-5^{31})\,x$,
using $5^{19}-2\frac{5^{19}-1}2=1$. We get similar formulas for $g_5$
and $g_6$ (replacing $19$ by resp. $25$ and $30$), etc.  Computing all
the $g_i$ this way, we get the following global triangular integral basis:
\[
  \{ 1,\frac \theta {5^4},\frac {\theta^2} {5^8},\frac {\theta^3}
  {5^{12}},\frac{g_4(\theta)}{5^{18}},\frac{g_5(\theta)}{5^{24}},\frac{g_6(\theta)}{5^{29}},\frac{g_7(\theta)}{5^{33}},\frac{g_8(\theta)}{5^{37}},\frac{g_9(\theta)}{5^{42}},\frac{g_{10}(\theta)}{2\,5^{51}},\frac{g_{11}(\theta)}{2\,5^{55}},\frac{g_{12}(\theta)}{2\,5^{59}}
  \}
\]
In particular, we get $D_f=2^3\,5^{372}$ mentioned earlier.

%%% Local Variables:
%%% mode: latex
%%% TeX-master: t
%%% End: